\documentclass[12pt]{article}

\def\Cov{\mathop{\rm Cov}}

\usepackage{amsmath,amsthm,amssymb,amsbsy,amsopn}
\usepackage{graphicx}

\newtheorem{Def}{Definition}[section]
\newtheorem{Thm}[Def]{Theorem}
\newtheorem{Lem}[Def]{Lemma}
\newtheorem{Prop}[Def]{Proposition}
\newtheorem{Coro}[Def]{Corollary}
\newtheorem{Rem}[Def]{Remark}
\theoremstyle{definition}

\newenvironment{enumerate2}{%
   \begin{list}%
   {%
   }%
   {%
      \usecounter{enum2}
      \setlength{\itemindent}{0em}
      \setlength{\leftmargin}{3em}
      \setlength{\rightmargin}{0em}
      \setlength{\labelsep}{1em}
      \setlength{\labelwidth}{3em}
      \setlength{\itemsep}{0em}
      \setlength{\parsep}{0em}
      \setlength{\listparindent}{0em}
   }
}{%
   \end{list}%
}

\title{On Fractional Tempered Stable Motion}
\author{C. Houdr\'e\footnote{Laboratoire d'Analyse et de Math\'ematiques
Appliqu\'ees, CNRS UMR 8050, Universit\'e Paris XII, 94010 Cr\'eteil
Cedex, France, and School of Mathematics, Georgia Institute of Technology,
Atlanta, GA, 30332-0160, USA, houdre@math.gatech.edu} 
~and R. Kawai\footnote{Quantitative Research Department, Daiwa Securities SMBC
Co.Ltd., 1-14-5, Eitai, Koto-ku, Tokyo, 135-0034, Japan,
reiichiro.kawai@daiwasmbc.co.jp}}
\date{February 28, 2005}
  
\begin{document}
\maketitle

\begin{abstract}
{\it Fractional tempered stable motion (fTSm)} is defined
 and studied.
FTSm has the same covariance structure as fractional Brownian motion,
 while having tails heavier than Gaussian but lighter than
 stable.
Moreover, in short time it is close to fractional stable
 L\'evy motion, while it is approximately fractional Brownian motion in long
 time.  
A series representation of fTSm is derived and used for simulation and
 to study some of its sample path properties.
\end{abstract}

\section{Introduction}

{\it Fractional Brownian motion (fBm)} and its various extensions are not only rich mathematical objects
but have also been extensively used in application to model asset price
dynamics, data traffic in telecommunication network, daily hydrological
series, and turbulence, to mention but a few topics.
We recall that standard fBm $\{B^H_t:t\in \mathbb{R}\}$, $H\in
(0,1]$ is a centered Gaussian process with continuous paths and with the
following covariance structure;
\begin{equation} 
\Cov (B^H_t,B^H_s)=\frac{1}{2}\left(|t|^{2H}+|s|^{2H}-|t-s|^{2H}\right), \quad
t,s\in\mathbb{R}.\label{covariance structure of fBm}
\end{equation}
However, modeling drawbacks of fBm and of some of its extensions
have also been discussed in the literature.
For example, although Gaussianity provides analytical
tractability, its light tails are often
inadequate for modeling higher variability observed in various natural
phenomena.
On the other hand, non-Gaussian stable generalizations
immediately lead to infinite second moment and to a lack of closed form
for the density, resulting in significant analytical difficulties.
Moreover, the selfsimilar and stationary increments properties of fBm are
sometimes unrealistic in practical modeling.

In order to remove these drawbacks, we introduce and study {\it
fractional tempered stable motion (fTSm)}, which has the following properties:
\begin{enumerate2}
\item[(i)] Its marginals have tails heavier than Gaussian but lighter than
(non-Gaussian) stable (Proposition \ref{ftsm marginal proposition}).
\item[(ii)] It has the same second order structure as fBm (Proposition
\ref{cov of fTSm}).
\item[(iii)] In long time it behaves like fBm, while in short time
	   it is more akin to fractional stable motions.
	   (Theorem \ref{fTSm to fSm}).
\end{enumerate2}
We present its series representation (Proposition \ref{fractional
Levy motion series representation}), which is potentially useful for
simulation and which is also used to study sample path properties.
When $H\in (1/\alpha ,1/\alpha +1/2)$, it has a.s. H\"older continuous sample
paths with exponent $(0,H-1/\alpha)$ (Proposition \ref{Holder continuity of
fTSm}) and is not a semimartingale (Proposition \ref{non-semimartingale}).
In contrast, the sample paths of fTSm become nowhere bounded as soon as
$H<1/\alpha$ (Proposition \ref{nowhere bounded}).

Let us close this section by introducing some notations and definitions
which will be used throughout the text.
$\mathbb{R}^d$ is the $d$-dimensional Euclidean space with the norm
$\|\cdot\|$.
$\mathbb{R}_0^d:=\mathbb{R}^d\setminus \{0\},$ and
$\mathcal{B}(\mathbb{R}^d_0)$ is the Borel $\sigma$-field of
$\mathbb{R}^d_0$.
$(\Omega ,\mathcal{F},\mathbb{P})$ is our underlying probability space.
$\mathcal{L}(Y)$ is the law of the random vector $Y$, while
``$\stackrel{\mathcal{L}}{=}$'' denotes equality in law, or
equality of the finite dimensional distributions when stochastic
processes are considered.
Similarly, ``$\stackrel{\mathcal{L}}{\to}$'' is used for convergence in
law, or of the finite dimensional distributions, while
``$\stackrel{d}{\to}$'' denotes the weak convergence of stochastic
processes in the space $D([0,\infty),\mathbb{R})$ of c\`adl\`ag
functions from $[0,\infty)$ into $\mathbb{R}$ equipped with the
Skorohod topology.
$C([0,\infty),\mathbb{R})$ is the space of continuous functions from
$[0,\infty)$ to $\mathbb{R}$ endowed with the uniform metric.
A sequence of stochastic processes $\{X^n_t:t\ge 0\}_{n\in
\mathbb{N}}$ in $C([0,\infty),\mathbb{R})$ is said to be
{\it tight} if for each compact set $K\in [0,\infty)$ and each $\epsilon
>0$,
\[
 \lim_{\delta \to 0}\limsup_{n\to\infty}\, \mathbb{P}\left(\sup_{t,s\in K,|t-s|\le \delta}|X^n_t-X^n_s|>\epsilon\right)=0.
\]
A sequence of stochastic processes $\{X^n_t:t\ge
 0\}_{n\in \mathbb{N}}$ is said to converge {\it uniformly on compacts in
 probability (ucp)} to a stochastic process $\{X_t:t\ge 0\}$ if for each
 compact set $K\in [0,\infty)$ and each $\epsilon >0$,
\[
\lim_{n\to \infty} \mathbb{P}\left(\sup_{t\in K}|X^n_t-X_t|>\epsilon\right)=0.
\]
This last convergence will be denoted by
``$X^n\stackrel{ucp}{\longrightarrow}X$''.
Finally, $\ln^+a=\ln a$ if $a\ge 1$.
$\ln^+a=0$ otherwise. 

\section{Definition of Fractional Tempered Stable Motion}\label{section def of fTSm}

We begin by briefly reviewing tempered stable distributions and
processes (Rosi\'nski~\cite{super rosinski}) which are building
blocks for the fTSm defined below.
Let $\mu$ be an infinitely divisible probability measure, without
Gaussian component, on $\mathbb{R}^d$.
Then, $\mu$ is called tempered stable if its L\'evy measure has the form
\[
 \nu
 (B)=\int_{\mathbb{R}_0^d}\int_0^{\infty}{\bf 1}_B(sx)s^{-\alpha-1}e^{-s}ds\rho
 (dx),\quad B\in \mathcal{B}(\mathbb{R}_0^d),
\]
where $\alpha \in (0,2)$ and where $\rho$, {\it the inner measure}, is
such that
\begin{equation}\label{original condition of inner measure}
 \int_{\mathbb{R}_0^d}\|x\|^{\alpha}\rho (dx)<\infty.
\end{equation}
The two parameters $\alpha$ and $\rho$ above uniquely identify the
L\'evy measure of tempered stable distributions.
Under the additional condition
\begin{equation*}
\begin{cases}
\int_{\mathbb{R}_0^d}\|x\|\rho (dx)<\infty,&{\rm if}~\alpha \in (0,1),\\
\int_{\mathbb{R}_0^d}\|x\|(1+\ln^+\|x\|)\rho (dx)<\infty,&{\rm if}~\alpha
=1,
\end{cases}
\end{equation*} 
the characteristic function of $\mu$ has a closed form expression given by
\begin{equation}\label{cf 1 of ts}
 \widehat{\mu}(y)=\exp\left[i\langle y,b \rangle+\int_{\mathbb{R}_0^d}
\phi_{\alpha}(\langle y,x\rangle )\rho (dx)\right],
\end{equation}
for some $b\in \mathbb{R}^d$ and where, for $s\in \mathbb{R},$
\begin{equation}\label{def of phi in ftsm}
 \phi_{\alpha}(s)=
\begin{cases}
\Gamma (-\alpha)((1-is)^{\alpha}-1+i\alpha s),&{\rm if}~ \alpha \in
 (0,1)\cup (1,2),\\
(1-is)\ln (1-is)+is,&{\rm if}~ \alpha =1.
\end{cases}
\end{equation}
Below, we write $\mu \sim TS (\alpha ,\rho ;b)$ if $\widehat{\mu}$
 is given by (\ref{cf 1 of ts}) and denote by $\{X^{TS}_t:t\ge 0\}$ a
 tempered stable L\'evy process in $\mathbb{R}$ such that
 $X^{TS}_1\sim TS(\alpha ,\rho ;b)$.
Setting $b=0$ gives $\mathbb{E}[X^{TS}_t]=0$ for every $t\ge
0$, and then $\{X^{TS}_t:t\ge 0\}$ is a martingale.

{\it From now on, we always assume that for any $t> 0$,
\[
\mathbb{E}[X^{TS}_t]=0,
\]
and further that
\begin{equation}
 \int_{\mathbb{R}_0}|x|^2\rho (dx)<\infty,\label{second moment condition}
\end{equation}
so that for any $t> 0$, $\mathbb{E}[(X^{TS}_t)^2]<+\infty$.}

We will define fTSm as a process of stochastic integral with respect to the
tempered stable process, i.e., $\{\int_{\mathcal{S}}f(t,s)dX^{TS}_s:t\in
\mathcal{T}\}$, where $f:\mathcal{T}\times \mathcal{S}\to \mathbb{R}$ is a
deterministic function.
Various such representations of fBm has been introduced in the
literature.
The moving-average representation (Mandelbrot and Van
Ness~\cite{mandelbrot van ness}) and the harmonizable representation 
are the most commonly used.
These have also been extended to non-Gaussian stable marginals.
(See Chapter 7 of Samorodnitsky and Taqqu~\cite{samorodnitsky taqqu}.)
A lesser known representation of fBm involves a Volterra kernel and is due to
Decreusefond and \"Ust\"unel~\cite{decreusefond ustunel}.
Recall that a Volterra kernel is a function $K:\mathbb{R}\times \mathbb{R}\to
\mathbb{R}$ such that $K(t,s)=0$ for $s>t$.
In the present paper, we will use the Volterra kernel
$K_{H,\alpha}:[0,\infty)\times [0,\infty)\to [0,\infty)$, given by
\begin{eqnarray}
K_{H,\alpha}(t,s)&:=&c_{H,\alpha}\bigg[\left(\frac{t}{s}\right)^{H-1/\alpha}(t-s)^{H-1/\alpha}\nonumber\\
&&\quad \qquad
 -\left(H-\frac{1}{\alpha}\right)s^{1/\alpha-H}\int_s^tu^{H-1/\alpha-1}(u-s)^{H-1/\alpha}du\bigg]\,{\bf 1}_{[0,t]}(s),\label{def of kernel}
\end{eqnarray}
where $H\in (1/\alpha-1/2,1/\alpha+1/2)$, $\alpha \in (0,2)$, and 
\[
 c_{H,\alpha}=\left(\frac{G(1-2G)\Gamma (1/2-G)}{\Gamma (2-2G)\Gamma (G+1/2)}\right)^{1/2},
\]
with (throughout) $G:=H-1/\alpha +1/2$.
Clearly, $K_{1/\alpha,\alpha}(t,s)= {\bf 1}_{[0,t]}(s)$.
Note that when $H \in (1/\alpha ,1/\alpha+1/2)$, we also have
\begin{equation}\label{easier representation}
 K_{H,\alpha}(t,s)=c_{H,\alpha}(H-1/\alpha)s^{1/\alpha
 -H}\int_s^t(u-s)^{H-1/\alpha -1}u^{H-1/\alpha}du\,{\bf 1}_{[0,t]}(s).
\end{equation}
Despite its complex structure, there are two main advantages to using
this kernel as an integrand:
\begin{enumerate2}
\item[(i)] It is defined only on $[0,t]$, $t>0$, while the domain of
	   definition of the moving-average kernel
	   $(t-s)_+^{H-1/\alpha}-(-s)_+^{H-1/\alpha}$ is the whole real
	   line.
(A moving-average fractional L\'evy motion has recently been studied
	   in Benassi {\it et al.} \cite{benassi}.)
           In general, however, it is impossible to generate background driving
	   stochastic processes defined on $\mathbb{R}$.

\item[(ii)] It can be written as a Riemann-Louiville
	   fractional integral, whose inverse function has a closed
	   form which is also a Riemann-Louiville fractional derivative.
This is important for the prediction of the sample paths of fTSm, problem of
	   consequence in financial modeling.
This will be presented in a subsequent paper \cite{finance paper}.
\end{enumerate2}

Below, we derive a necessary and sufficient condition
on $H$ so that for each $t>0$ and $p\ge 2$, the kernel is in $L^p([0,t])$.

\begin{Lem}\label{Volterra kernel integrability}
Let $t>0$, let $\alpha \in (0,2)$, and let $p\ge 2$.
Then, $K_{H,\alpha}(t,\cdot)\in L^p([0,t])$ if and only if $H\in
 (1/\alpha-1/p,1/\alpha+1/p).$
In particular, $K_{H,\alpha}(t,\cdot)\in L^2([0,t])$.
Moreover, when $K_{H,\alpha}(t,\cdot)\in L^p([0,t])$,
\begin{equation}
 \int_0^tK_{H,\alpha}(t,s)^pds=C_{H,\alpha,p}\,t^{p(H-1/\alpha)+1},\label{def of kernel const}
\end{equation}
where 
\begin{equation}\label{kernel L^p constant}
 C_{H,\alpha,p}=c_{H,\alpha}^p\int_0^1v^{p\left(\frac{1}{\alpha} -H\right)}\left[(1-v)^{H-\frac{1}{\alpha}}-\left(H-\frac{1}{\alpha}\right)\int_v^1
w^{H-\frac{1}{\alpha}-1}(w-v)^{H-\frac{1}{\alpha}}dw\right]^pdv.
\end{equation}
\end{Lem}

\begin{proof}
The case $H=1/\alpha$ is trivial since $K_{H,\alpha}(t,s)={\bf 1}_{[0,t]}(s)$.
If $H>1/\alpha$, then $K_{H,\alpha}(t,s)\ge 0$, $K_{H,\alpha}(t,s)$ is
 decreasing in $s$, and $K_{H,\alpha}(t,s)\sim C's^{1/\alpha -H}$ as $s
 \downarrow 0$ for some constant $C'$. 
Hence, $K_{H,\alpha}(t,\cdot)\in L^p([0,t])$ if and only if $p(1/\alpha
 -H)>-1$, i.e., $H<1/\alpha +1/p$.
When $H<1/\alpha$, $K_{H,\alpha}(\cdot,s)$ explodes at $s=0$ and $s=t$.
In fact, $K_{H,\alpha}(t,s)\sim C''s^{H-1/\alpha}$ as
 $s\downarrow 0$ and $K_{H,\alpha}(t,s)\sim C'''(t-s)^{H-1/\alpha}$ as
 $s\uparrow t$ for some constants $C''$ and $C'''$.
Thus, $K_{H,\alpha}(t,\cdot)\in L^p([0,t])$ if and only if
 $p(H-1/\alpha)>-1$, i.e., $H>1/\alpha-1/p$, which proves the first claim.
The last claim follows from elementary computation.
\end{proof}

\begin{Rem}\label{instead p}{\rm
Above, we only considered the case $p\ge 2$ because of the moment
condition (\ref{second moment condition}) and since $H\in (1/\alpha
-1/2,1/\alpha +1/2)$ in (\ref{def of kernel}).
However, it is easily seen that the results above remain true for
 arbitrary $p>0$, provided that the kernel is defined for
 arbitrary $H$.
In particular, we have $K_{H,\alpha}(t,\cdot)\in L^{\alpha}([0,t])$
 since $(0,2/\alpha)\supset (1/\alpha-1/2,1/\alpha+1/2).$ 
}\end{Rem}

Let us state two other known properties of the kernel.
The proof of (ii) below can be found in, e.g., Decreusefond and \"Ust\"unel
\cite{decreusefond ustunel}, or Nualart \cite{nualart 2}, while (i) is
immediate.

\begin{Lem}\label{Volterra selfsimilar}
(i) For each $h>0$,
\[
 K_{H,\alpha}(ht,s)=h^{H-1/\alpha}K_{H,\alpha}(t,s/h).
\]
(ii) For $t,s>0$,
\[
 \int_0^{t\land s}K_{H,\alpha}(t,u)K_{H,\alpha}(s,u)du=\frac{1}{2}(t^{2G}+s^{2G}-|t-s|^{2G}),
\]
where $G=H-1/\alpha+1/2$.
\end{Lem}

We are now in a position to define fTSm.

\begin{Def}
Fractional tempered stable motion $\{L_t^H:t\ge 0\}$ in $\mathbb{R}$ is
 given by
\begin{equation}
L_t^H:=\int_0^tK_{H,\alpha}(t,s)dX^{TS}_s,\quad t\ge 0,\label{def of fTSm} 
\end{equation} 
where the integral is defined in the $L^2(\Omega
 ,\mathcal{F},\mathbb{P})$-sense.
\end{Def}

The integral above is well defined by the moment condition
(\ref{second moment condition}), Lemma \ref{Volterra
selfsimilar} (ii) and the help of the Wiener-It\^o isometry.
(See also the proof of Proposition \ref{cov of fTSm} below.)
For convenience, we will henceforth write $\{L^H_t:t\ge 0\}\sim fTSm
(H,\alpha ,\rho)$ when $\{L^H_t:t\ge 0\}$ is defined as (\ref{def of fTSm}).
We note that $L^{1/\alpha}_{\cdot}=X^{TS}_{\cdot}$, which is a L\'evy
process, since $K_{1/\alpha,\alpha}(t,s)= {\bf 1}_{[0,t]}(s)$.

\vspace{1em}
The following is an important result on the marginals of fTSm.

\begin{Prop}\label{ftsm marginal proposition}
The finite dimensional distributions of fTSm are tempered stable with
 finite second moment.
\end{Prop}

\begin{proof}
Let $k\in \mathbb{N}.$
It suffices to show that for any real sequence
 $\{a_i\}_{i=1}^k$ and any nonnegative nondecreasing real sequence
 $\{t_i\}_{i=1}^k$, the random variable $\sum_{i=1}^ka_iL^H_{t_i}$ is
 tempered stable.
First, observe that
\[
 \sum_{i=1}^ka_iL^H_{t_i}=\int_0^{t_k}\left(\sum_{i=1}^ka_iK_{H,\alpha}(t_i,s)\right)dX^{TS}_s.
\]
Then, by Proposition 35 of Rocha-Arteaga and Sato~\cite{rocha-arteaga sato},
 we get
\begin{eqnarray*}
\mathbb{E}[e^{iy\sum_{i=1}^ka_iL^H_{t_i}}]&=&\exp\left[\int_0^{t_k}\int_{\mathbb{R}_0}\phi_{\alpha}(yx\sum_{i=1}^ka_iK_{H,\alpha}(t_i,s))\rho(dx)ds\right]\\
&=&\exp\left[\int_{\mathbb{R}_0}\phi_{\alpha}(yx)\eta(dx)\right],
\end{eqnarray*}
where $\phi_{\alpha}$ is given by (\ref{def of phi in ftsm}) and where $\eta =M\circ J$ with $M(dx,ds)=\rho (dx)ds$ and
\[
  J(B)=\left\{(x,s)\in \mathbb{R}_0\times [0,t_k]:x\sum_{i=1}^{k}a_i
K_{H,\alpha} (t_i,s)\in B\right\}, \quad B\in\mathcal{B}(\mathbb{R}_0).
\]
The measure $\eta$ is well defined as an inner measure with finite
 second moment since for each $i$, $K_{H,\alpha}(t_i,\cdot)\in
 L^2([0,t_i])$ and
\[
\int_{\mathbb{R}_0}|x|^2\eta(dx)=\int_{\mathbb{R}_0}
|x|^2\rho (dx)\,\int_0^{t_k}\left(\sum_{i=1}^ka_iK_{H,\alpha}(t_i,s)
\right)^2ds<\infty,
\]
which concludes the proof.
\end{proof}

It is worth noting the one dimensional marginal result as a corollary.

\begin{Coro}\label{one dimensional marginal}
Let $\{L^H_t:t\ge 0\}\sim fTSm (H,\alpha,\rho)$ and let $\phi_{\alpha}$
 be given by (\ref{def of phi in ftsm}).
For each $t>0$,
\begin{equation}\label{cf of fTSm}
 \mathbb{E}[e^{iyL^H_t}]=\exp \left[\int_0^t\int_{\mathbb{R}_0}\phi_{\alpha}(yxK_{H,\alpha}(t,s))\rho(dx)ds\right],
\end{equation}
and thus
\begin{equation}\label{ftsm marginal is ts}
 L^H_t\sim TS(\alpha, \eta_t;0)
\end{equation}
where $\eta_t =M\circ J_t$ with $M(dx,ds)=\rho (dx)ds$ and $J_t(B)=\{(x,s)\in
 \mathbb{R}_0\times [0,t]:xK_{H,\alpha}(t,s)\in B\}$, $B\in\mathcal{B}(\mathbb{R}_0)$.
\end{Coro}

\section{Covariance Structure and Long-range Dependence}

Let us first describe the covariance structure of fTSm.
 
\begin{Prop}\label{cov of fTSm}
Let $\{L^H_t:t\ge 0\}\sim fTSm(H,\alpha ,\rho)$. Then, $\mathbb{E}[L^H_t]=0,$ and
\begin{equation}
\Cov
 (L_t^H,L_s^H)=\frac{1}{2}\left(t^{2G}+s^{2G}-|t-s|^{2G}\right)\mathbb{E}[(X^{TS}_1)^2],\quad
 s,t>0.\label{covariance structure of fTSm}
\end{equation}
\end{Prop}

\begin{proof}
Recall that $\{X^{TS}_t:t\ge 0\}$ is a square-integrable centered
 martingale.
The first claim is thus trivial.
For the second claim, observe that for $s\in [0,t]$,
\begin{eqnarray*}
 \Cov
  (L^H_t,L^H_s)&=&\mathbb{E}[L^H_tL^H_s]\\
&=&\mathbb{E}\left[\int_0^tK_{H,\alpha}(t,u)dX^{TS}_u\int_0^sK_{H,\alpha}(s,u)dX^{TS}_u\right]\\
&=&\mathbb{E}\left[\int_0^{s\land
	      t}K_{H,\alpha}(t,u)dX^{TS}_u\int_0^{s\land t}K_{H,\alpha}(s,u)dX^{TS}_u\right]\\
&=&\mathbb{E}[(X^{TS}_1)^2]\,\int_0^{s\land t}K_{H,\alpha}(t,u)K_{H,\alpha}(s,u)du,
\end{eqnarray*}
where the last equality holds by the Wiener-It\^o isometry.
Lemma \ref{Volterra selfsimilar} (ii) then gives the result.
\end{proof} 

Let us state some immediate consequences of the previous result.

\begin{Coro}\label{second-order stuff corollary}
Let $\{L^H_t:t\ge 0\}\sim fTSm (H,\alpha,\rho)$.
For each $t> 0$ and each $h> 0,$
\begin{equation}\label{second-order selfsimilarity}
 \mathbb{E}[(L^H_{ht})^2]=h^{2G}\mathbb{E}[(L^H_t)^2],
\end{equation}
and for $s,t>0$,
\begin{equation}\label{second-order stationary increments}
 \mathbb{E}[(L^H_t-L^H_s)^2]=\mathbb{E}[(L^H_{|t-s|})^2]=|t-s|^{2G}\mathbb{E}[(X^{TS}_1)^2].
\end{equation}
\end{Coro}

\vspace{1em}
The property (\ref{second-order selfsimilarity}) is sometimes
 called {\it second-order selfsimilarity}.
Moreover, (\ref{second-order stationary increments}) says that fTSm
 has {\it second-order stationary increments}, which clearly
implies its continuity in probability.

We are now in a position to discuss the long-range dependence of fTSm.
The definition of long-range dependence is often ambiguous and varies
among authors.
In the present paper, we will follow Samorodnitsky and Taqqu
\cite{samorodnitsky taqqu}; the increments of a second-order stochastic
process $\{X_t:t\ge 0\}$ exhibit {\it long-range dependence} if for each
$h>0$,
\[
\sum_{n=1}^{\infty}|\Cov(X_h-X_0,X_{nh}-X_{(n-1)h})|=\infty,
\]
or {\it short-range dependence}, if each $h>0$, 
\[
\sum_{n=1}^{\infty}|\Cov(X_h-X_0,X_{nh}-X_{(n-1)h})|<\infty.
\]

\begin{Prop}\label{long-range dependence}
The increments of fTSm exhibit long-range dependence
 when $H\in (1/\alpha,1/\alpha+1/2)$, and short-range dependence when $H\in
 (1/\alpha -1/2,1/\alpha]$.
\end{Prop}

\begin{proof}
By Lemma \ref{cov of fTSm}, we have for each $h>0$,
\begin{eqnarray*}
\Cov(L_h^H,L_{t+h}^H-L_t^H)&=&\frac{1}{2}t^{2G}((1+h/t)^{2G}-2+(1-h/t)^{2G})\\
&\sim& \frac{1}{2}t^{2(G-1)} G(2G-1)h^2,
\end{eqnarray*}
as $t\to \infty$.
The claim then holds since $2(G-1)>-1$ for $H\in (1/\alpha,1/\alpha +1/2)$,
 while $2(G-1)\le -1$ for $H\in (1/\alpha-1/2,1/\alpha]$.
\end{proof}

In relation to the second moment, we will consider higher moments of fTSm.

\begin{Prop}\label{higher moments of fTSm}
Let $\{L^H_t:t\ge 0\}\sim fTSm(H,\alpha,\rho).$
Then, for each $p>2$ and each $t> 0$, $\mathbb{E}[|L^H_t|^p]<\infty$ if
 and only if $H\in (1/\alpha -1/p,1/\alpha +1/p)$ and
 $\int_{|x|>1}|x|^p\rho(dx)<\infty$.
\end{Prop}

\begin{proof}
By Corollary \ref{one dimensional marginal}, for each $t\ge 0$,
 $\mathcal{L}(L^H_t)$ is tempered stable.
By Proposition 2.3 (iii) of Rosi\'nski~\cite{super rosinski},
 $\mathbb{E}[|L^H_t|^p]<\infty$ if and only if $\int_{|x|>1}|x|^p\eta_t
 (dx)<\infty,$ where $\eta_t$ is the inner measure of $L^H_t$ given as in
 Corollary \ref{one dimensional marginal}, that is,
\[
 \iint_{|x|K_{H,\alpha}(t,s)>1}(|x|K_{H,\alpha}(t,s))^p\rho (dx)ds<\infty.
\]
The left hand side of the above can be decomposed into two terms;
\[
\int_0^tK_{H,\alpha}(t,s)^p\int_{\frac{1}{K_{H,\alpha}(t,s)}<|x|\le
 \frac{1}{K_{H,\alpha}(t,s)}\lor 1}|x|^p\rho (dx)ds,
\]
and
\[
 \int_0^tK_{H,\alpha}(t,s)^p\int_{|x|>\frac{1}{K_{H,\alpha}(t,s)}\lor
 1}|x|^p\rho(dx)ds.
\]
The first term is equivalent to $\int_0^tK_{H,\alpha}(t,s)^pds$ due to
 the moment condition (\ref{second moment condition}) on $\rho$, while
 the second terms is clearly equivalent to
 $\int_0^tK_{H,\alpha}(t,s)^pds \int_{|x|>1}|x|^p\rho (dx).$
Then, Lemma \ref{Volterra kernel integrability} concludes the proof.
\end{proof}

\begin{Rem}{\rm
It is shown in Proposition 2.3 (iv) of Rosi\'nski \cite{super rosinski}
that a tempered stable distribution has exponential moment of
certain order if and only if its inner measure has a compact support.
Unfortunately, the tempered stable marginal of fTSm cannot have
 exponential moment since, with the notation of the preceding
 proposition, for any $\epsilon >0$,
\[
 \eta_t(\{x\in
 \mathbb{R}_0:|x|>\epsilon\})=\int_{|x|K_{H,\alpha}(t,s)>\epsilon}\rho
 (dx)ds>0,
\]
due to the unboundedness of the kernel $K_{H,\alpha}$.
}\end{Rem}

\section{Series Representation}

In this section, we derive a series representation of fTSm, which is
inherited from the one of tempered stable processes obtained by Rosi\'nski
\cite{super rosinski}.
This representation can also be used for simulation. (See Figure \ref{fTSm and TS via series representation}.)
Moreover, we will make use of its structure to derive some sample path
properties in Section \ref{path property section}.

Let $\{\Gamma_i\}_{i\ge 1}$ be arrival times of a standard Poisson process,
let $\{E_i\}_{i\ge 1}$ be a sequence of iid exponential random variables
with parameter $1$, let $\{U_i\}_{i\ge 1}$ be a sequence of iid uniform
random variables on $[0,1]$, let $\{V_i\}_{i\ge 1}$ be a sequence of iid
random variables in $\mathbb{R}_0$ with common distribution 
\[
 \frac{|x|^{\alpha}\rho (dx)}{m(\rho)^{\alpha}},
\]
and let $\{T_i\}_{i\ge 1}$ be a sequence of iid uniform random variables on
$[0,T]$.
Also let $m(\rho)^{\alpha}$, $k'$, and $z_T$ be constants given by
$m(\rho)^{\alpha}=\int_{\mathbb{R}_0}|x|^{\alpha}\rho (dx)$, $k'=m(\rho)^{-\alpha}\int_{\mathbb{R}_0}x|x|^{\alpha-1}\rho(dx)$, and
\[
 z_T=
\begin{cases}
m(\rho)(\alpha/T)^{-1/\alpha}\zeta(1/\alpha)k'T^{-1}+|\Gamma
 (1-\alpha)|\int_{\mathbb{R}_0}x\rho (dx),&{\rm if}~\alpha \ne 1,\\
(\ln (m(\rho)T)+2\gamma)\int_{\mathbb{R}_0}x\rho
 (dx)-\int_{\mathbb{R}_0}x\ln|x|\rho(dx),&{\rm if}~\alpha =1,
\end{cases}
\]
where $\zeta$ denotes the Riemann zeta function and $\gamma
(=0.5772...)$ is the Euler constant.
Then, Theorem 5.4 of Rosi\'nski~\cite{super rosinski} tells us that
\[
\sum_{i=1}^{\infty}\Bigg[\left(m(\rho )
\left(\frac{\alpha\Gamma_i}{T}\right)^{-1/\alpha}\land E_iU_i^{1/\alpha}
|V_i|\right)\frac{V_i}{|V_i|}\,{\bf 1}(T_i\le t)-m(\rho)\left(\frac{\alpha
\,i}{T}\right)^{-1/\alpha}k'\frac{t}{T}\Bigg]+z_T\,t
\]
converges a.s. uniformly in $t\in [0,T]$ to a tempered stable process with
$TS(\alpha ,\rho ;0).$
This series representation can easily be extended to fTSm as follows.

\begin{Prop}\label{fractional Levy motion series representation}
Let $\{L^H_t:t\ge 0\}\sim fTSm(H,\alpha ,\rho)$ and let $T>0$.
Then,
\begin{eqnarray}
\nonumber &&\{L^H_t:t\in [0,T]\}\\
\nonumber &&\quad \stackrel{\mathcal{L}}{=}\Bigg\{\sum_{i=1}^{\infty}
\Bigg[\left(m(\rho )\left(\frac{\alpha\Gamma_i}{T}\right)^{-1/\alpha}
\land E_iU_i^{1/\alpha}|V_i|\right)\frac{V_i}{|V_i|}K_{H,\alpha}(t,T_i)\\
&&\qquad \quad -m(\rho)\left(\frac{\alpha \,i}{T}\right)^{-1/\alpha}k'
C_{H,\alpha,1}\frac{t^{H-1/\alpha+1}}{T}\Bigg]+z_TC_{H,\alpha,1}t^{H-1/\alpha
+1}:t\in [0,T]
\Bigg\},\label{series representaion for fTSm}
\end{eqnarray}
where $C_{H,\alpha,1}$ is the constant defined by (\ref{kernel L^p constant}).
If $\rho$ is symmetric, then 
\[
\{L^H_t:t\in [0,T]\}\stackrel{\mathcal{L}}{=}\Bigg\{\sum_{i=1}^{\infty}
\left(m(\rho )\left(\frac{\alpha\Gamma_i}{T}\right)^{-1/\alpha}
\land E_iU_i^{1/\alpha}|V_i|\right)\frac{V_i}{|V_i|}K_{H,\alpha}(t,T_i):t\in [0,T]
\Bigg\}.
\]
Moreover, if $H\in [1/\alpha,1/\alpha+1/2),$ then the series converges
 almost surely uniformly in $t$ to $fTSm(H,\alpha,\rho)$.
\end{Prop}

\begin{proof}
We will only consider the asymmetric case.
By arguments as in Theorem 5.4 of Rosi\'nski~\cite{super rosinski},
we get
\[
 \sum_{i=1}^{\infty}\left[m(\rho)\left(\frac{\alpha \,i}{T}
\right)^{-1/\alpha}k'C_{H,\alpha,1}\frac{t^{H-1/\alpha
+1}}{T}-c_i(T)\mathbb{E}[K_{H,\alpha}(t,T_1)]\right]=z_TC_{H,\alpha,1}t^{H-1/\alpha +1},
\] 
uniformly in $t$, where 
\[
 c_i(T):=\int_{i-1}^i\mathbb{E}\left[\left(m(\rho )
\left(\frac{\alpha r}{T}\right)^{-1/\alpha}\land E_1U_1^{1/\alpha}
|V_1|\right)\frac{V_1}{|V_1|}\right]dr.
\]
Hence, the right hand side of (\ref{series representaion for fTSm}) can be
 rewritten as
\begin{equation}\label{ftsm series Z}
\sum_{i=1}^{\infty}\left[\left(m(\rho )
\left(\frac{\alpha\Gamma_i}{T}\right)^{-1/\alpha}\land E_iU_i^{1/\alpha}
|V_i|\right)\frac{V_i}{|V_i|}K_{H,\alpha}(t,T_i)-c_i(T)\mathbb{E}[K_{H,\alpha}
(t,T_1)]\right]:=Z_t.
\end{equation}
Next, we need to analyze the finite dimensional distributions of
 $\{Z_t:t\in [0,T]\}.$
Let $k\in \mathbb{N},$ let $\{a_j\}_{j=1}^k$ be a real sequence, and let
 $\{t_j\}_{j=1}^k$ be a nondecreasing sequence taking values in $[0,T]$.
We will show that the random variable $\sum_{j=1}^ka_jZ_{t_j}$ has
 the same law as $\sum_{j=1}^ka_jL^H_{t_j}$.
In view of Proposition \ref{ftsm marginal proposition}, we have
\[
 \mathbb{E}\left[e^{iy\sum_{j=1}^ka_jL^H_{t_j}}\right]=\exp\left[\int_0^T
\int_{\mathbb{R}_0}\phi_{\alpha}(yx\sum_{j=1}^ka_jK_{H,\alpha}(t_j,s))\rho
 (dx)ds\right],
\]
where $\phi_{\alpha}$ is given by (\ref{def of phi in ftsm}).
Also, observe that
\begin{eqnarray*}
&& \sum_{j=1}^ka_jZ_{t_j}=\sum_{i=1}^{\infty}\Bigg[\left(m(\rho )
\left(\frac{\alpha\Gamma_i}{T}\right)^{-1/\alpha}\land E_iU_i^{1/\alpha}
|V_i|\right)\frac{V_i}{|V_i|}\sum_{j=1}^ka_jK_{H,\alpha}(t_j,T_i)\\
&&\qquad \qquad \qquad \qquad \qquad \qquad \qquad \qquad \qquad 
-c_i(T)\mathbb{E}\Bigg[\sum_{j=1}^ka_jK_{H,\alpha}(t_j,T_1)\Bigg]\Bigg].
\end{eqnarray*}
This series representation is induced by the L\'evy measure 
\begin{eqnarray*}
&&\nu (B)\\
&&=\int_0^T\int_{\mathbb{R}_0}\int_0^{\infty}\int_0^1\int_0^{\infty}{\bf
 1}_{B}(H(r/T,u,s,x)\sum_{j=1}^ka_jK_{H,\alpha}(t_j,v))drdue^{-s}ds
\rho_1(dx)\frac{dv}{T}\\
&&=\int_0^T\int_{\mathbb{R}_0}\int_0^{\infty}{\bf 1}_B(sx\sum_{j=1}^ka_j
K_{H,\alpha}(t_j,v))s^{-\alpha -1}e^{-s}ds\rho (dx)dv,
\end{eqnarray*}
where $H(r,u,s,x)=(m(\rho)(\alpha r)^{-1/\alpha}\land
 su^{1/\alpha}|x|) x/|x|$ and
 $\rho_1(dx)=m(\rho)^{-\alpha}|x|^{\alpha}\rho (dx).$
In fact, the measure $\nu$ is well defined as a L\'evy measure since
 $K_{H,\alpha}(t_j,\cdot)\in L^2([0,t_j])$ for each $j$.
Therefore, by Theorem 4.1 (B) of Rosi\'nski~\cite{rosinski2},
\begin{eqnarray*}
 \mathbb{E}[e^{iy\sum_{j=1}^ka_jZ_{t_j}}]&=&\exp\left[\int_{\mathbb{R}_0}(e^{iyz}-1-iyz)\nu(dz)\right]\\
&=&\exp\left[\int_0^T
\int_{\mathbb{R}_0}\phi_{\alpha}(yx\sum_{j=1}^ka_jK_{H,\alpha}(t_j,s))\rho
 (dx)ds\right],
\end{eqnarray*}
which proves the equality of all finite dimensional distributions.

Next, let $H \in [1/\alpha ,1/\alpha+1/2)$. 
Define, for $t\in [0,T]$ and $s\in [0,\infty)$,
\begin{eqnarray*}
&& Z_{t,s}:=\sum_{\{i:\Gamma_i\le s\}}\Bigg[\left(m(\rho )
\left(\frac{\alpha\Gamma_i}{T}\right)^{-1/\alpha}\land E_iU_i^{1/\alpha}
|V_i|\right)\frac{V_i}{|V_i|}K_{H,\alpha}(t,T_i)\\
&&\qquad \qquad \qquad \qquad \qquad \qquad \qquad \qquad \qquad \qquad 
-c_i(T)\mathbb{E}[K_{H,\alpha}(t,T_1)]\Bigg].
\end{eqnarray*}
Notice that $\{Z_{t,s}:t\in [0,T],\,s\in [0,\infty)\}$ has independent
 increments in $s$ ({\it not} in $t$, of course.)
Arguments as in the proof of Theorem 5.1 of \cite{rosinski2} give the
 a.s. convergence of the series uniformly on $[0,T]$. 
\end{proof}

\section{Sample Path Properties}\label{path property section}

In this section, we investigate sample path properties of fTSm.
Let us begin with the case $H\in (1/\alpha-1/2,1/\alpha)$.

\begin{Prop}\label{nowhere bounded}
When $H\in (1/\alpha-1/2,1/\alpha)$, fTSm is a.s. unbounded on every
 interval of positive length.
\end{Prop}

\begin{proof}
Let $T>0$.
For each $t\in (0,T]$, $\lim_{s\uparrow t}K_{H,\alpha}(t,s)=+\infty$.
Hence, in the series representation given in Proposition \ref{fractional
 Levy motion series representation}, $K_{H,\alpha}(T_i,T_i)=+\infty$, for
 all $i\in \mathbb{N}$, and so none of the summands are well defined.
This implies that $\sup_{t\in [0,T]}|L^H_t|=\infty$ a.s.
\end{proof}

Unfortunately, the above sample path property makes fTSm with
short-range dependence of little practical use.
We notice that for any $H\in (1/\alpha-1/2,1/\alpha+1/2)$,
$K_{H,\alpha}(t,0)=+\infty$.
But, this turns out to be irrelevant to the unboundedness of sample
paths since $T_i\ne 0$ a.s., $i\in\mathbb{N}$.

FTSm with long-range dependence has better sample path properties.
In particular, it has a H\"older continuous version with exponent $\gamma \in
(0,H-1/\alpha)$.

\begin{Prop}\label{Holder continuity of fTSm}
If $H\in (1/\alpha ,1/\alpha+1/2),$ there exists a continuous
 modification of fTSm, which
 is a.s. locally H\"older continuous with exponent $\gamma$ for every
 $\gamma \in (0,H-1/\alpha)$.
\end{Prop}

\begin{proof}
By Corollary \ref{second-order stuff corollary}, we have $\mathbb{E}[|L^H_t-L^H_s|^2]=|t-s|^{2G}\mathbb{E}[(X^{TS}_1)^2]$.
If $H>1/\alpha,$ then $2G>1$, and thus the Kolmogorov-\v{C}entsov
 Theorem (see, for example, Theorem 3.23 of
 Kallenberg~\cite{kallenberg}) directly applies, giving the result.
\end{proof}

{\it We will henceforth always assume that when $H\in (1/\alpha ,1/\alpha
+1/2)$, we are using such a continuous version of fTSm.}

\begin{Prop}\label{non-semimartingale}
Let $\{L^H_t:t\ge 0\}\sim fTSm (H,\alpha,\rho)$ with $H\in (1/\alpha ,1/\alpha+1/2).$

(i) 
\[
 \lim_{N\to\infty}\mathbb{E}\left[\sum_{n=0}^{N-1}|L^H_{\frac{n+1}{N}T}-L^H_{\frac{n}{N}T}|^2\right]=
 0.
\]

(ii) FTSm is a.s. of infinite variation on every interval of positive length.

(iii) FTSm is not semimartingale.
\end{Prop}

\begin{proof}

(i) Immediate from
\[
 \mathbb{E}\left[\sum_{n=0}^{N-1}|L^H_{\frac{n+1}{N}T}-L^H_{\frac{n}{N}T}|^2\right]=(N/T)^{
2(1/\alpha -H)}.
\]

(ii) Let $T>0$. For each $s\in [0,T]$, we have
\[
 \limsup_{t_1,t_2\downarrow s}\frac{|K_{H,\alpha}(t_1,s)-K_{H,\alpha}(t_2,s)|}{|t_1-t_2|}=+\infty,
\]
which implies that for each $s\in [0,T]$, $K_{H,\alpha}(\cdot,s)$ is of
 infinite variation.
By Theorem 4 of Rosi\'nski~\cite{rosinski path} and a symmetrization
 argument given there, fTSm is of infinite variation with
 positive probability.
FTSm is selfdecomposable and hence by Corollary 3 of Rosi\'nski
 \cite{rosinski zero-one law}, it obeys a zero-one law.
This gives the result.

(iii) 
The convergence in (i) implies convergence in probability, and together
 with (ii), the claim follows from the same argument in Lin \cite{lin}.
\end{proof}

\begin{Rem}{\rm
In view of (iii) above, the stochastic integration for fTSm cannot be
 defined in the classical semimartingale sense.
However, a slight modification of the kernel $K_{H,\alpha}$ induces 
a corresponding semimartingale, which can be arbitrarily close to fTSm. 
Recall the simpler expression of the kernel for $H\in
 (1/\alpha,1/\alpha+1/2)$ given by (\ref{easier representation}).
Observe that
\[
 \frac{\partial}{\partial
 t}K_{H,\alpha}(t,s)=c_{H,\alpha}(H-1/\alpha)(t-s)^{H-1/\alpha-1}
\left(\frac{t}{s}\right)^{H-1/\alpha}{\bf 1}_{[0,t)}(s),
\]
and that
\[
 K_{H,\alpha}(t,s)=\int_s^t\frac{\partial}{\partial
 u}K_{H,\alpha}(u,s)du\,{\bf 1}_{[0,t)}(s).
\]
Therefore,
\begin{equation}
 L^H_t=\int_0^tK_{H,\alpha}(t,s)dX^{TS}_s=\int_0^t\left(\int_s^t
\frac{\partial}{\partial
u}K_{H,\alpha}(u,s)du\right)dX^{TS}_s.\label{fTSm double integral form}
\end{equation}
If the two integrals could be interchanged, then fTSm would be
of finite variation, i.e., it would be a semimartingale; we have seen
 that this is not the case.
On the other hand, the integrability condition of the stochastic
Fubini's theorem (Theorem 46 of Protter~\cite{protter book})
can be achieved by slightly modifying the kernel $K_{H,\alpha}$.
Set 
\begin{equation}\label{regularized volterra kernel}
 K_{H,\alpha}^n(t,s):=c_{H,\alpha}(H-1/\alpha)s^{1/\alpha
  -H}\int_s^t\left(u+\frac{1}{n}-s\right)^{H-1/\alpha-1}u^{H-1/\alpha}du\,{\bf 1}_{[0,t]}(s),\quad n\in \mathbb{N},
\end{equation}
and so we have
\[
 \frac{\partial}{\partial
 t}K^n_{H,\alpha}(t,s)=c_{H,\alpha}(H-1/\alpha)\left(t+\frac{1}{n}-s\right)^{H-1/\alpha-1}\left(\frac{t}{s}\right)^{H-1/\alpha}{\bf
 1}_{[0,t]}(s).
\]
The integrability condition is then satisfied; for every $u\in [0,t]$,
\begin{eqnarray*}
&&\int_0^u\left(\frac{\partial}{\partial u}K^n_{H,\alpha}(u,s)\right)^2ds\\
&&\qquad \qquad =c_{H,\alpha}^2(H-1/\alpha)^2\int_0^u\left(u+\frac{1}{n}-s\right)^{2(H-1/\alpha -1)}\left(\frac{u}{s}\right)^{2(H-1/\alpha)}ds\\
&&\qquad \qquad \le c_{H,\alpha}^2(H-1/\alpha)^2(1-2(H-1/\alpha))^{-1}n^{-2(H-1/\alpha -1)}u <\infty ,
\end{eqnarray*}
and thus the stochastic Fubini's theorem applies.
Therefore, we get
\[
 \int_0^t
 K_{H,\alpha}^n(t,s)dX^{TS}_s=\int_0^t\left(\int_0^u\frac{\partial}{\partial
 u}K^n_{H,\alpha}(u,s)dX^{TS}_s\right)du,
\]
which is clearly of finite variation.
For financial modeling, it is of interest to further modify the
 above to get an infinite variation semimartingale.
This can be done as follows.
For $\epsilon >0$, set $K^{n,\epsilon}_{H,\alpha}(t,s):=K^n_{H,\alpha}(t,s)+\epsilon.$
Since $\frac{\partial}{\partial t}K^{n,\epsilon}_{H,\alpha}(t,s)=
\frac{\partial}{\partial t}K^n_{H,\alpha}(t,s),$ the stochastic Fubini's
 theorem again applies and thus we get 
\begin{eqnarray*}
 \int_0^t K_{H,\alpha}^{n,\epsilon}(t,s)dX^{TS}_s&=&\int_0^t(\epsilon +K^n_{H,\alpha}(t,s))dX^{TS}_s\\
&=&\epsilon X^{TS}_t+\int_0^t\left(\int_0^u\frac{\partial}{\partial u}K^n_{H,\alpha}(u,s)dX^{TS}_s\right)du,
\end{eqnarray*}
which exactly follows the definition of the canonical decomposition of
semimartingales, i.e. a martingale plus a finite variation process.
}\end{Rem}

\section{Short and Long Time Behavior}

In this section, we obtain the {\it short and
long time behavior} of fTSm, which are also inherited from the background
driving tempered stable processes obtained by Rosi\'nski \cite{super
rosinski}.
In short time, fTSm is asymptotically fractional stable motion (to be
defined below), while in long time it is approximately fBm.

Let us begin by briefly reviewing the corresponding behaviors of
tempered stable processes, proved in Theorem 3.1 of \cite{super rosinski}.

\vspace{1em}
(i) \underline{{\it Short time behavior}}:
Let $\{X^{TS}_t:t\ge 0\}\sim TS(\alpha,\rho ;0)$ and let
\begin{equation}\label{def of b for fsm}
 b_{h,\alpha}=
\begin{cases}
h\Gamma(1-\alpha)\int_{\mathbb{R}_0}x\rho (dx),&if~\alpha \in (0,1),\\
-(1+\ln h)\int_{\mathbb{R}_0}x\rho (dx),&if~\alpha =1,\\
0,&if~\alpha \in (1,2).
\end{cases}
\end{equation} 
Then,
\[
 \left\{h^{-1/\alpha}\left(X^{TS}_{ht}+b_{h,\alpha}t\right):t\ge
 0\right\}\stackrel{d}{\to}\{X^{\alpha}_t:t\ge 0\},\quad as~h\to 0,
\]
where $\{X^{\alpha}_t:t\ge 0\}$ is an $\alpha$-stable process in
$\mathbb{R}$ such that
\[
 \mathbb{E}[e^{iyX^{\alpha}_t}]=\exp\left[t\int_{\mathbb{R}_0}\varphi_{\alpha}(yx)\rho
 (dx)\right],
\]
where
\begin{equation}\label{def of varphi for fsm}
\varphi_{\alpha}(s)=
\begin{cases}
-\Gamma(-\alpha)\cos \frac{\pi\alpha}{2}|s|^{\alpha}(1-i\tan
 \frac{\pi\alpha}{2}{\rm sgn} (s)),&{\rm if}~\alpha \in (0,1)\cup (1,2),\\
-\left(\frac{\pi}{2}|s|+is\ln |s|\right)+is,&{\rm if}~\alpha =1.
\end{cases} 
\end{equation}

\vspace{1em}
(ii) \underline{{\it Long time behavior}}: 
Let $\{X^{TS}_t:t\ge 0\}\sim TS(\alpha,\rho ;0)$.
Then,
\[
 \{h^{-1/2}X_{ht}^{TS}:t\ge 0\}\stackrel{d}{\to}\{cW_t:t\ge 0\},\quad
 as~h\to \infty,
\]
where $\{W_t:t\ge 0\}$ is a standard (centered) Brownian motion and
\begin{equation}\label{def of c for fBm}
 c^2=\Gamma (2-\alpha)\int_{\mathbb{R}_0}x^2\rho (dx).
\end{equation}

We will say that the limiting stable process $\{X^{\alpha}_t:t\ge 0\}$
given above is {\it associated to} the tempered stable process
$\{X^{TS}_t:t\ge 0\}$.

Let us now define {\it fractional stable motions (fSm),} which turns out
to be a short time limit of fTSm.
Below, for $\alpha \ne 1$, the integral in (\ref{def of fSm}) is well
defined in probability since $K_{H,\alpha}(t,\cdot)\in
L^{\alpha}([0,t]),$ $t>0$.
On the other hand, when $\alpha =1$, the extra symmetry assumption on
$\rho$ also ensures that it is well defined in probability.
(See Remark \ref{instead p} and Samorodnitsky and Taqqu \cite
{samorodnitsky taqqu}.)

\begin{Def}
Let $\{L^H_t:t\ge 0\}\sim fTSm(H,\alpha,\rho)$, where when $\alpha =1$,
 $\rho$ is additionally assumed to be symmetric.
Fractional stable motion (fSm) $\{L_t^{H,\alpha}:t\ge 0\}$
 associated to fTSm $\{L^H_t:t\ge 0\}$ is given via 
\begin{equation}
L_t^{H,\alpha}:=\int_0^tK_{H,\alpha}(t,s)dX^{\alpha}_s,\quad t\ge 0.
\label{def of fSm}
\end{equation}
\end{Def}

Let us derive some basic properties of fSm.

\begin{Lem}\label{fSm marginal}
Let $\{L^{H,\alpha}_t:t\ge 0\}$ be fSm associated to $\{L^H_t:t\ge
 0\}\sim fTSm(H,\alpha,\rho).$

(i) The finite dimensional distributions of fSm are stable.

(ii) For each $t>0$, the characteristic function of $L^{H,\alpha}_t$ is
 given by
\begin{equation}\label{fsm cf}
\mathbb{E}[e^{iyL^{H,\alpha}_t}]=\exp\left[C_{H,\alpha,\alpha}t^{\alpha H}
\int_{\mathbb{R}_0}\widetilde{\varphi}_{\alpha}(yx)\rho (dx)\right],
\end{equation}
where $C_{H,\alpha,\alpha}$ is the constant given by (\ref{kernel L^p
 constant}) and where
\[
 \widetilde{\varphi}_{\alpha}(s)=
\begin{cases}
\varphi_{\alpha}(s),&if~\alpha \in (0,1)\cup (1,2),\\
-\frac{\pi}{2}|s|,&if~\alpha =1,
\end{cases}
\]
where $\varphi_{\alpha}$ is defined by (\ref{def of varphi for fsm}).

(iii) FSm is selfsimilar; $\{h^{-H}L^{H,\alpha}_{ht}:t\ge 0\}
\stackrel{\mathcal{L}}{=}\{L^{H,\alpha}_t:t\ge 0\}.$

(iv) FSm has (strictly) stationary increments;
for each $t>s>0$, $L^{H,\alpha}_t-L^{H,\alpha}_s\stackrel{\mathcal{L}}{=}
L^{H,\alpha}_{t-s}.$

(v) When $H\in (1/\alpha,1/\alpha +1/2)$, fSm has a continuous
 version, while when $H\in (1/\alpha-1/2,1/\alpha)$, it is unbounded on
 every interval of positive length.

(vi) When $H=1/\alpha$, fSm is an $\alpha$-stable (L\'evy) process.

(vii) With the notation in Section 4,
 if $\alpha \in (1,2)$, then
\begin{eqnarray*}
&&\{L^{H,\alpha}_t:t\in [0,T]\}\\
&&\quad \stackrel{\mathcal{L}}{=}
\Bigg\{\sum_{i=1}^{\infty}\Bigg[m(\rho)\left(\frac{\alpha
 \Gamma_i}{T}\right)^{-1/\alpha}\frac{V_i}{|V_i|}K_{H,\alpha}(t,T_i)-m(\rho)\left(\frac{\alpha i}{T}\right)^{-1/\alpha}k'C_{H,\alpha,1}\frac{t^{H-1/\alpha+1}}{T}\Bigg]\\
&&\quad \qquad \qquad \qquad \qquad \qquad \qquad +m(\rho)\left(\frac{\alpha}{T}\right)^{-1/\alpha}\zeta
 (1/\alpha)k'C_{H,\alpha,1}\frac{t^{H-1/\alpha+1}}{T}:t\in [0,T]\Bigg\},
\end{eqnarray*}
while if $\alpha \in (0,1],$ or if $\alpha \in (1,2)$ and $\rho$ is
 symmetric, then
\[
\{L^{H,\alpha}_t:t\in [0,T]\}\stackrel{\mathcal{L}}{=}
\Bigg\{\sum_{i=1}^{\infty}m(\rho)\left(\frac{\alpha
 \Gamma_i}{T}\right)^{-1/\alpha}\frac{V_i}{|V_i|}K_{H,\alpha}(t,T_i):
t\in [0,T]\Bigg\}.
\]
Moreover, if $H\in [1/\alpha,1/\alpha+1/2)$, $\alpha \in (0,2),$ the
 above series converges almost surely uniformly in $t\in [0,T]$ to a
 version of $\{L^{H,\alpha}_t:t\in [0,T]\}$. 
\end{Lem}

\begin{proof}
With the notation of Proposition \ref{ftsm marginal proposition},
for $\alpha \ne 1$,
\begin{eqnarray*}
&&\mathbb{E}[e^{iy\sum_{i=1}^ka_iL^{H,\alpha}_{t_i}}]\\
&&=\exp\Bigg[-\Gamma(-\alpha)\cos\frac{\pi \alpha}{2}\int_0^{t_k}
\int_{\mathbb{R}_0}\left|\sum_{i=1}^ka_iK_{H,\alpha}(t_i,s)
yx\right|^{\alpha}\\
&&\qquad \qquad \qquad \qquad  \left(1-i\tan
\frac{\pi\alpha}{2}{\rm
sgn}\left(\left|\sum_{i=1}^ka_iK_{H,\alpha}(t_i,s)
\right|^{\alpha}yx\right)\right)\rho (dx)ds\Bigg]\\
&&=\exp\Bigg[-\Gamma(-\alpha)\cos\frac{\pi \alpha}{2}\int_0^{t_k}
\left|\sum_{i=1}^ka_iK_{H,\alpha}(t_i,s)\right|^{\alpha}ds\\
&&\qquad \qquad \qquad \qquad \qquad \qquad \qquad 
 \int_{\mathbb{R}_0}|yx|^{\alpha}\left(1-i\tan \frac{\pi\alpha}{2}
{\rm sgn}(yx)\right)\rho (dx)\Bigg].
\end{eqnarray*}
When $\alpha =1,$ the symmetry of $\rho$ yields $b_{h,1}=0,$ and so
\begin{eqnarray*}
\mathbb{E}[e^{iy\sum_{i=1}^ka_iL^{H,1}_{t_i}}]&=&\exp\left[-\int_0^t\int_{\mathbb{R}_0}
\left(\frac{\pi}{2}\left|\sum_{i=1}^ka_iK_{H,1}(t_i,s)yx\right|\right)\rho(dx)ds\right]\\
&=&\exp\Bigg[-\frac{\pi}{2}\int_0^t\left|\sum_{i=1}^ka_iK_{H,1}(t_i,s)
\right|ds\,\int_{\mathbb{R}_0}|yx|\rho(dx)\Bigg].
\end{eqnarray*}
For each $\alpha \in (0,2)$, $\sum_{i=1}^ka_iK_{H,\alpha}(t_i,\cdot)\in
 L^{\alpha}([0,t_k])$ since $a_iK_{H,\alpha}(t_i,\cdot)\in
 L^{\alpha}([0,t_i])$ for each $i.$
Thus, (i) holds.  
Clearly, (ii) is a direct consequence of (i).
(iii) and (iv) follow from the selfsimilarity and stationary increments
 properties of $\{X^{(\alpha)}_t:t\ge 0\}$ with the help of Lemma
 \ref{Volterra selfsimilar}.
By (iii) and (iv),
 $E[|L^{H,\alpha}_t-L^{H,\alpha}_s|^p]=|t-s|^{pH}E[|L^{H,\alpha}_1|^p]$,
 $0<p<\alpha$.
Hence, when $H\in (1/\alpha,1/\alpha+1/2)$, the continuity follows from the
 Kolmogorov-\v{C}entsov Theorem.
Next, let $H\in (1/\alpha -1/2,1/\alpha)$.
For each $T>0$, $\sup_{t\in [0,T]}|K_{H,\alpha}(t,s)|=+\infty$, $s\in [0,T].$
The nowhere boundedness follows from (i) and Theorem 4 of
 Rosi\'nski~\cite{rosinski path} as well as the zero-one law for stable
 processes and a symmetrization argument given there.
Hence, (v) holds.
(vi) is immediate from $K_{1/\alpha,\alpha}(t,s)={\bf 1}_{[0,t]}(s).$
Finally, (vii) follows from arguments as in Proposition 5.5 of
 Rosi\'nski \cite{super rosinski} and from Proposition \ref{fractional
 Levy motion series representation}. 
\end{proof}

{\it We will henceforth always assume that when $H\in
(1/\alpha ,1/\alpha +1/2)$, we are using a continuous version of fSm.}

\begin{Rem}{\rm
Many extensions of fBm are available in the stable literature;
for example, linear fractional stable motion, log-fractional stable
 motion, harmonizable fractional stable motion.
(See, e.g., Samorodnitsky and Taqqu \cite{samorodnitsky taqqu}.)
These various extensions are not necessarily identical in law since their
 marginals are determined by kernels in the stochastic integral
 representations.
Indeed, fSm defined above is still different from any of them.
}\end{Rem}

We are now in a position to present the main result of this section.

\begin{Thm}\label{fTSm to fSm}
Let $\{L^H_t:t\ge 0\}\sim fTSm (H,\alpha,\rho)$ with $H\ne 1/\alpha$.

(i) \underline{Short time behavior}: Let 
\[
 b=
\begin{cases}
\Gamma (1-\alpha) \int_{\mathbb{R}_0}x\rho (dx),&if~\alpha \in
 (0,1),\\
0,&if~\alpha \in [1,2),
\end{cases}
\]
and let $k_t=\int_0^tK_{H,\alpha}(t,s)ds$.
Then,
\[
 \{h^{-H}L^H_{ht}+h^{1-1/\alpha}bk_t:t\ge
 0\}\stackrel{\mathcal{L}}{\to}\{L^{H,\alpha}_t:t\ge 0\}\quad as~h\to 0,
\]
where $\{L^{H,\alpha}_t:t\ge
 0\}$ is fSm associated to $\{L^H_t:t\ge 0\}.$

\vspace{1em}
(ii) \underline{Long time behavior}:
\[
 \{h^{-G}L^H_{ht}:t\ge 0\}\stackrel{\mathcal{L}}{\to}\{c B^{G}_t:t\ge
 0\}\quad as ~ h\to \infty,
\]
where $\{B^G_t:t\ge 0\}$ is a standard fBm with parameter $G$ and where
 $c$ is the constant given by (\ref{def of c for fBm}).

(iii) When $H\in (1/\alpha ,1/\alpha +1/2)$, the convergence in (i) and
 (ii) can be strengthened to the weak convergence in
 $C([0,\infty),\mathbb{R}).$  
\end{Thm}

\begin{proof}[Proof of (i) and (ii)]
(i) Observe that for each $t\ge 0$,
\[
h^{-H}L^H_{ht}+h^{1-1/\alpha}b\,k_t
=\int_0^{t}K_{H,\alpha}(t,s)h^{-1/\alpha}d(X^{TS}_{hs}+b\,hs):=Y_t^h.
\]
It thus suffices to show that for any real sequence $\{a_i\}_{i=1}^k$
 and nonnegative nondecreasing real sequence $\{t_i\}_{i=1}^k$,
 $k\in \mathbb{N},$ the random variable $\sum_{i=1}^ka_iY_{t_i}^h$
 converges in law to $\sum_{i=1}^ka_iL^{H,\alpha}_{t_i}$, as $h\to 0$.
Since
\[
 \sum_{i=1}^ka_iY_{t_i}^h=\int_0^{t_k}\left(\sum_{i=1}^ka_iK_{H,\alpha}
(t_i,s)\right) h^{-1/\alpha}d(X^{TS}_{hs}+b\,hs),
\]
we have by Proposition~\ref{ftsm marginal proposition} that
\begin{equation}
 \mathbb{E}\left[e^{iy\sum_{i=1}^ka_iY_{t_i}^h}\right]=\exp\left[\int_0^{t_k}\int_{\mathbb{R}_0}h\psi_{\alpha}(yxh^{-1/\alpha}\sum_{i=1}^ka_iK_{H,\alpha}(t_i,s))\rho
 (dx)ds\right],\label{before taking limit for fSm}
\end{equation}
where 
\begin{equation}\label{inner function}
 \psi_{\alpha}(s)=
\begin{cases} 
\Gamma (-\alpha)((1-is)^{\alpha}-1),&{\rm if}~\alpha \in (0,1),\\
\frac{1}{2}\ln (1+s^2)-s\tan^{-1}s,&{\rm if}~\alpha =1,\\
\Gamma (-\alpha)((1-is)^{\alpha}-1+i\alpha s),&{\rm if}~\alpha \in (1,2).
\end{cases}
\end{equation}
Note that $\psi_1$ is obtained via the symmetry of $\rho$.
(See Proposition 2.8 of Rosi\'nski \cite{super rosinski}.)
We then want to show that (\ref{before taking limit for fSm})
 tends to the characteristic function of the random variable $\sum_{i=1}^ka_iL^{H,\alpha}_{t_i}$, as $h\to 0.$

The proof of Theorem 3.1 (i) of Rosi\'nski~\cite{super rosinski} shows
 that for $\alpha \ne 1$,
\[
 \lim_{h\to 0}h \psi_{\alpha}(h^{-1/\alpha}s)=\varphi_{\alpha}(s),
\]
where $\varphi_{\alpha}$ is given by (\ref{def of varphi for fsm}) and  
\[
 \left|h\psi_{\alpha}(h^{-1/\alpha}s)\right|\le z_{\alpha}|s|^{\alpha},
\]
where $z_{\alpha}$ is some constant depending only on $\alpha$($\ne 1$).
When $\alpha =1$, we have 
\[
 \lim_{h\to 0}h\psi_{1}(h^{-1}s)=-\frac{\pi}{2}|s|,
\]
and the uniform boundedness (in $h>0$) of $|h\psi_1(h^{-1}s)|$ can be
 shown as
\begin{eqnarray*}
|h\psi_1(h^{-1}s)|&\le& |h\ln \sqrt{1+h^{-2}s^2}|+|s \tan^{-1}(h^{-1}s)|\\
&\le&|h\ln (1+h^{-1}|s|)|+\frac{\pi}{2}|s|\\
&\le & \left(1+\frac{\pi}{2}\right)|s|.
\end{eqnarray*}
Clearly, $\sum_{i=1}^ka_iK_{H,\alpha}(t_i,s)\in L^{\alpha}([0,t_k])$
 since $a_iK_{H,\alpha}(t_i,s)$ are in $L^{\alpha}([0,t_i])$. 
Together with the moment condition (\ref{original condition of inner
 measure}) on $\rho$, the passage to the limit in
 (\ref{before taking limit for fSm}) is justified.
Hence,
\[
 \lim_{h\to
 0}\mathbb{E}\left[e^{iy\sum_{i=1}^ka_iY_{t_i}^h}\right]=\exp\left[\int_0^{t_k}\int_{\mathbb{R}_0}\varphi_{\alpha}(yx\sum_{i=1}^ka_iK_{H,\alpha}(t_i,s))\rho
 (dx)ds\right],
\]
which is the characteristic function of $\sum_{i=1}^ka_iL^{H,\alpha}_{t_i}.$

(ii) We have that for each $h>0$, 
\[
 \Cov
 (h^{-G}L^H_{ht},h^{-G}L^H_{hs})=\frac{1}{2}\left(t^{2G}+s^{2G}-(t-s)^{2G}\right)\mathbb{E}[(X_1^{TS})^2],\quad
 s\in [0,t].
\]
Hence, for the convergence of all finite dimensional distributions, we
 only need to show that the marginal law at any time of
 $\{h^{-G}L^H_{ht}:t\ge 0\}$ converges to Gaussian.
Without loss of generality, let $t=1$.
By Lemma \ref{Volterra selfsimilar} (i),
\begin{eqnarray*}
 \mathbb{E}[e^{iyh^{-G}L^H_h}]&=&\exp\left[\int_0^h\int_{\mathbb{R}_0}
\vartheta_{\alpha} (h^{-G}yxK_{H,\alpha}(h,s))\rho(dx)ds\right]\\
&=&\exp\left[\int_0^1\int_{\mathbb{R}_0}\vartheta_{\alpha}
(h^{-1/2}yxK_{H,\alpha}(1,s))\rho(dx)ds\right],
\end{eqnarray*}
where $\vartheta_{\alpha}(u)=\int_0^{\infty}(e^{ius}-1-ius)s^{-\alpha-1}
e^{-s}ds.$
As in the proof of Theorem 3.1 (ii) of Rosi\'nski~\cite{super
 rosinski}, it follows that
\[
 |\vartheta_{\alpha}(yxK_{H,\alpha}(1,s))|\le (yx)^2\Gamma(2-\alpha)\int_0^1K_{H,\alpha}(1,s)^2ds,
\] 
which justifies the passage to the limit below
\[
 \lim_{h\to
 \infty}\mathbb{E}[e^{iyh^{-G}L^H_h}]=\exp\left[-\frac{y^2}{2}\int_0^1K_{H,\alpha}(1,s)^2ds\Gamma
 (2-\alpha)\int_{\mathbb{R}_0}x^2\rho (dx)\right].
\]
This shows the convergence to a Gaussian law, which concludes the proof.
\end{proof}

To proof (iii), let us present a technical lemma.

\begin{Lem}\label{tightness lemma}
Let $\{X^n_t:t\ge 0\}_{n\in \mathbb{N}}$ and $\{Y^n_t:t\ge 0\}_{n\in
 \mathbb{N}}$ be sequences of stochastic processes in
 $C([0,\infty),\mathbb{R})$.
If the sequence $\{X^n_t:t\ge 0\}_{n\in \mathbb{N}}$ is tight and if for
 each $n\in\mathbb{N}$, there exists a version $\widetilde{Y}^n$of $Y^n$
 defined on the same probability space as $X^{n}$ and such that the sequence
 $\{\widetilde{Y}^n_t-X^n_t:t\ge 0\}$ converges ucp to
 zero, then the sequence $\{Y^n_t:t\ge 0\}_{n\in \mathbb{N}}$ is tight.
\end{Lem}

\begin{proof}
For each compact set $K\in [0,\infty)$ and for each $\delta >0,$ we have
\begin{eqnarray}
&& \nonumber \mathbb{E}\left[\sup_{t,s\in K,|t-s|\le
	      \delta}\left|Y^{n}_t-Y^{n}_s
\right|\land 1\right]\\
&&\nonumber \quad \le \mathbb{E}\left[\sup_{t,s\in K,|t-s|\le \delta}\left|\widetilde{Y}^{n}_t-X^{n}_t\right|\land
 1\right]+\mathbb{E}\left[\sup_{t,s\in K,|t-s|\le
 \delta}\left|X^{n}_t-X^{n}_s\right|\land 1\right]\\
&&\nonumber \qquad \qquad +\mathbb{E}\left[\sup_{t,s\in K,|t-s|\le \delta}
\left|X^{n}_s-\widetilde{Y}^{n}_s\right|\land 1\right]\\
&&\quad = 2\,\mathbb{E}\Bigg[\sup_{t\in K}\left|\widetilde{Y}^{n}_t-
X^{n}_t\right|\land 1\Bigg]+\mathbb{E}\left[\sup_{t,s\in K,|t-s|\le
 \delta}\left|X^{n}_t-X^{n}_s\right|\land 1\right].\label{tightness
lemma equation}
\end{eqnarray}
The first term in (\ref{tightness lemma equation}) tends to zero as
 $n\to \infty$, since
 $\widetilde{Y}^{n}-X^{n}\stackrel{ucp}{\longrightarrow} 0$.
The second term also tends to zero as
 $n\to \infty$ and $\delta \to 0$ by the tightness of $\{X^n_t:t\ge
 0\}_{n\in\mathbb{N}}$ in $C([0,\infty),\mathbb{R}).$
The claimed result then follows.
\end{proof}

\begin{proof}[Proof of Theorem \ref{fTSm to fSm} (iii)]
By, for example, Lemma 16.2, Theorem 16.3 and 16.5 of Kallenberg
 \cite{kallenberg}, it suffices to show the tightness of the sequences
 $\{h^{-H}L^H_{ht}+h^{1-1/\alpha}bk_t:t\ge 0\}$ (as $h\downarrow 0$) and
 $\{h^{-G}L^H_{ht}:t\ge 0\}$ (as $h\uparrow \infty$) in
 $C([0,\infty),\mathbb{R})$.

We begin with the short time behavior case.
By Lemma \ref{fSm marginal} (iii) and (iv), for each $p\in (0,\alpha),$
\begin{equation}\label{tightness of fSm}
 E[|h^{-H}L^{H,\alpha}_{ht}-h^{-H}L^{H,\alpha}_{hs}|^p]=(t-s)^{pH}E[|L^{H,\alpha}_1|],\quad
 s\in [0,t].
\end{equation}
By Corollary 16.9 of Kallenberg \cite{kallenberg}, the uniform
 boundedness in $h$ seen in (\ref{tightness of fSm}) implies the
 tightness of the sequence $\{h^{-H}L^{H,\alpha}_{ht}:t\ge 0\}$ in
 $C([0,\infty),\mathbb{R})$.
Hence, by Lemma \ref{tightness lemma}, it is enough to find
 $h^{-H}L^{H,\alpha}_{h\cdot}$ and
 $h^{-H}L^H_{h\cdot}+h^{1-1/\alpha}bk_{\cdot}$ in $C([0,T],\mathbb{R}),$
 $T>0$, defined on a common probability space and such that
 $h^{-H}L^{H,\alpha}_{h\cdot}-(h^{-H}L^H_{h\cdot}+h^{1-1/\alpha}bk_{\cdot})
\stackrel{ucp}{\longrightarrow}0$.
To this end, we make use of their series representations.
First, let $\alpha \in (0,1]$, or let $\alpha \in (1,2)$ with a
 symmetric $\rho$.
With the notation of Section 4, it follows from Proposition
 \ref{fractional Levy motion series representation} and Lemma \ref{fSm
 marginal} (vii) that the stochastic processes
\[
h^{-H}\sum_{i=1}^{\infty}m(\rho)\left(\frac{\alpha
 \Gamma_i}{T}\right)^{-1/\alpha}\frac{V_i}{|V_i|}K_{H,\alpha}(h\cdot,hT_i):=h^{-H}\widetilde{L}^{H,\alpha}_{h\cdot}
\]
and 
\[
h^{-H}\sum_{i=1}^{\infty}\left(m(\rho)\left(\frac{\alpha
 \Gamma_i}{T}\right)^{-1/\alpha}\land E_iU_i^{1/\alpha}|V_i|\right)\frac{V_i}{|V_i|}K_{H,\alpha}(h\cdot,hT_i):=h^{-H}\widetilde{L}^H_{h\cdot}+h^{1-1/\alpha}b k_{\cdot},
\]
converge almost surely uniformly on $[0,T]$, respectively, to versions of
 fSm $L^H_{\cdot}$ and fTSm $h^{-H}L^{H,\alpha}_{h\cdot}+h^{1-1/\alpha}b
 k_{\cdot}$, defined on a common probability space by using the common
 random sequences $\{\Gamma_i\}_{i\ge 1},$ $\{T_i\}_{i\ge 1},$
 $\{V_i\}_{i\ge 1},$ $\{E_i\}_{i\ge 1}$ and $\{U_i\}_{i\ge 1}$.
Then, in view of Lemma \ref{Volterra selfsimilar} (i),
\begin{eqnarray}
&& \nonumber h^{-H}\widetilde{L}^{H,\alpha}_{h\cdot}-
(h^{-H}\widetilde{L}^H_{h\cdot}+h^{1-1/\alpha}b k_{\cdot})\\
&&\qquad =\sum_{i=1}^{\infty}\left[\left(m(\rho)\left(\frac{\alpha
 \Gamma_i}{T}\right)^{-1/\alpha}-h^{-1/\alpha}E_iU_i^{1/\alpha}|V_i|\right)\lor 0\right]\frac{V_i}{|V_i|}K_{H,\alpha}(\cdot,T_i),\label{difference of fSm and fTSm}
\end{eqnarray}
which clearly converges ucp to zero as $h\to 0$.
Also, in case $\alpha \in (1,2)$ and $\rho$ is asymmetric, a
 similar argument yields (\ref{difference of fSm and fTSm}).
Hence, the sequence $\{h^{-H}L^H_{ht}+h^{1-1/\alpha}bk_t:t\ge 0\}$ is
 tight in $C([0,\infty),\mathbb{R})$, which concludes the short time
 behavior case.

For the long time behavior case, Corollary \ref{second-order stuff
 corollary} gives for $h>0,$
\begin{equation}\label{tightness}
 \mathbb{E}[(h^{-G}L^H_{ht}-h^{-G}L^H_{hs})^2]=(t-s)^{2G}\mathbb{E}[(X_1^{TS})^2].
\end{equation}
Again by Corollary 16.9 of Kallenberg \cite{kallenberg},
the uniform boundedness in $h$ seen in (\ref{tightness})
 implies the tightness of the sequence $\{h^{-G}L^H_{ht}:t\ge 0\}$ in
 $C([0,\infty),\mathbb{R})$, which completes the proof.
\end{proof}

\begin{Rem}{\rm
The short time behavior result can also be seen from the series representation.
For simplicity, consider the symmetric case.
With the notations of Theorem \ref{fractional Levy motion series
 representation}, we have by Lemma \ref{Volterra selfsimilar} (i)
\begin{eqnarray}
 h^{-H}L^H_{ht}&\stackrel{\mathcal{L}}{=}&
h^{-H}\sum_{i=1}^{\infty}\left(m(\rho )
\left(\frac{\alpha\Gamma_i}{hT}\right)^{-1/\alpha}\land E_iU_i^{1/\alpha}
|V_i|\right)\frac{V_i}{|V_i|}K_{H,\alpha}(ht,hT_i)\nonumber \\
&=&\sum_{i=1}^{\infty}\left(m(\rho )
\left(\frac{\alpha\Gamma_i}{T}\right)^{-1/\alpha}\land E_iU_i^{1/\alpha}
|V_i|h^{-1/\alpha}\right)\frac{V_i}{|V_i|}K_{H,\alpha}(t,T_i)\nonumber \\
&\to&\sum_{i=1}^{\infty}m(\rho )\left(
\frac{\alpha\Gamma_i}{T}\right)^{-1/\alpha}\frac{V_i}{|V_i|}
K_{H,\alpha}(t,T_i)\quad a.s.,\quad {\rm as}~h\to 0.
\end{eqnarray}
}\end{Rem}

\vspace{1em}
We have seen in Lemma \ref{fSm marginal} that the marginals of fSm
are $\alpha$-stable and so their covariance does not exist.
An alternative notion is the one of {\it covariation}.
For two jointly symmetric
$\alpha$-stable random variables $X$ and $Y$ with $\alpha >1$;
\begin{equation}\label{taqqu def of codifference}
\tau (X,Y):=\|X\|^{\alpha}_{\alpha}+\|Y\|^{\alpha}_{\alpha}
-\|X-Y\|^{\alpha}_{\alpha},
\end{equation} 
where the norm $\|\cdot \|$ gives the scale of parameter, i.e. for
$Z\sim S_{\alpha}(\sigma ,0,0)$, $\|Z\|_{\alpha}=\sigma$.
More generally, one can also define the {\it codifference} for any
jointly infinitely divisible random variables $X$ and $Y$ as
\begin{equation}
I(\theta_1 ,\theta_2;X,Y):=-\ln E[e^{i(\theta_1 X+\theta_2 Y)}]+\ln
E[e^{i\theta_1 X}]+\ln E[e^{i\theta_2 Y}],\label{def of rosinski codifference}
\end{equation}
for $\theta_1,\theta_2\in \mathbb{R}.$
Clearly, (\ref{taqqu def of codifference}) is a special case of (\ref{def of
rosinski codifference}).
Let $\{L^{H,\alpha}_t:t\ge 0\}$ be fSm associated to $\{L^H_t:t\ge
0\}\sim fTSm (H,\alpha,\rho)$ where $\rho$ is {\it symmetric}.
Then,
\begin{equation}
I (1,-1;L^{H,\alpha}_t,L^{H,\alpha}_s)=C(t^{\alpha H}+s^{\alpha
 H}-(t-s)^{\alpha H}), \quad s\in [0,t],
\end{equation}
for some constant $C$.
In the Gaussian case, the codifference coincides with the covariance. 
For example, for a standard fBm $\{B^G_t:t\ge 0\}$,
\[
\tau (B^G_t,B^G_s)=\frac{1}{2}(t^{2G}+s^{2G}-(t-s)^{2G})=\Cov
(B^G_t,B^G_s),\quad s\in [0,t],
\]
and
\[
\tau (B^G_{t+1}-B^G_t,B^G_1-B^G_0)=\Cov (B^G_{t+1}-B^G_t,B^G_1-B^G_0)\sim
C t^{2(G-1)},
\]
as $t\to\infty$.
See Samorodnitsky and Taqqu \cite{samorodnitsky taqqu} for more details of the covariation and the codifference.

Interestingly enough, as shown below, the codifference of the increments of fTSm has the
same order of decay as the covariance (see Proposition \ref{long-range dependence}).

\begin{Prop}\label{codifference of fTSm}
Let $\{L^H_t:t\ge 0\}\sim fTSm(H,\alpha ,\rho)$.
Then, 
\[
I(\theta_1,\theta_2;L^H_{t+1}-L^H_t,L^H_1-L^H_0)\sim
C(\theta_1,\theta_2) t^{2(G-1)},
\]
as $t\to \infty$, where 
\[
C(\theta_1,\theta_2)=\frac{-ic_{H,\alpha}\theta_1\pi}{\Gamma (\alpha )\sin (\pi \alpha )}\int_{[0,1]\times \mathbb{R}_0}\left((1-ix\theta_2K_H(1,s))^{\alpha -1}-1\right)xs^{1/\alpha -H}ds \rho (dx).
\]
\end{Prop}

\begin{proof}
Observe that
\begin{eqnarray*}
&&I(\theta_1,\theta_2;L^H_{t+1}-L^H_t,L^H_1-L^H_0)\\
&&\qquad
=\Gamma (-\alpha )\int_{\mathbb{R}\times\mathbb{R}_0}\Big(-(1-ix(\theta_1(K_{H,\alpha}(t+1,s)-K_{H,\alpha}(t,s))+\theta_2K_{H,\alpha}(1,s)))^{\alpha}\\
&&\qquad \qquad \qquad \qquad \qquad +(1-ix\theta_1(K_{H,\alpha}(t+1,s)-K_{H,\alpha}(t,s)))^{\alpha}\\
&&\qquad \qquad \qquad \qquad \qquad \qquad \qquad \qquad \qquad
+(1-ix\theta_2K_{H,\alpha}(1,s))^{\alpha}-1\Big)ds\rho (dx),
\end{eqnarray*}
and that
\[
 K_{H,\alpha}(t+1,s)-K_{H,\alpha}(t,s)\sim c_{H,\alpha} s^{1/\alpha
 -H}t^{2(G-1)},
\] 
as $t\to \infty$.
Hence, for each $s>0$,
\begin{eqnarray*}
&&-(1-ix(\theta_1(K_{H,\alpha}(t+1,s)-K_{H,\alpha}(t,s))+\theta_2K_{H,\alpha}(1,s)))^{\alpha}\\
&&\qquad \qquad +(1-ix\theta_1(K_{H,\alpha}(t+1,s)-K_{H,\alpha}(t,s)))^{\alpha}+(1-ix\theta_2K_{H,\alpha}(1,s))^{\alpha}-1\\
&&\quad \sim
i\alpha x\theta_1c_{H,\alpha}s^{1/\alpha -H}\left((1-ix\theta_2K_{H,\alpha}(1,s))^{\alpha-1}-1\right)t^{2(G-1)},
\end{eqnarray*}
as $t\to\infty.$
The result then holds since $\Gamma
 (-\alpha)=\frac{-\pi}{\alpha \Gamma (\alpha )\sin (\pi \alpha )}.$ 
\end{proof}

\section{Concluding Remarks}

Willinger, Taqqu and Teverovsky \cite{willinger taqqu teverovsky} assert
that a numerical analysis of stock price time series indicates long-range
dependence with marginal tails heavier than Gaussian but
lighter than stable.
Moreover, it has been known that in shorter time, the asset
price paths tend to lack higher moments, while they have a
Gaussian behavior in long time.
Indeed, fTSm achieves all those properties.

\begin{figure}
\begin{center}
\begin{tabular}{cc}
      \resizebox{65mm}{!}{\includegraphics{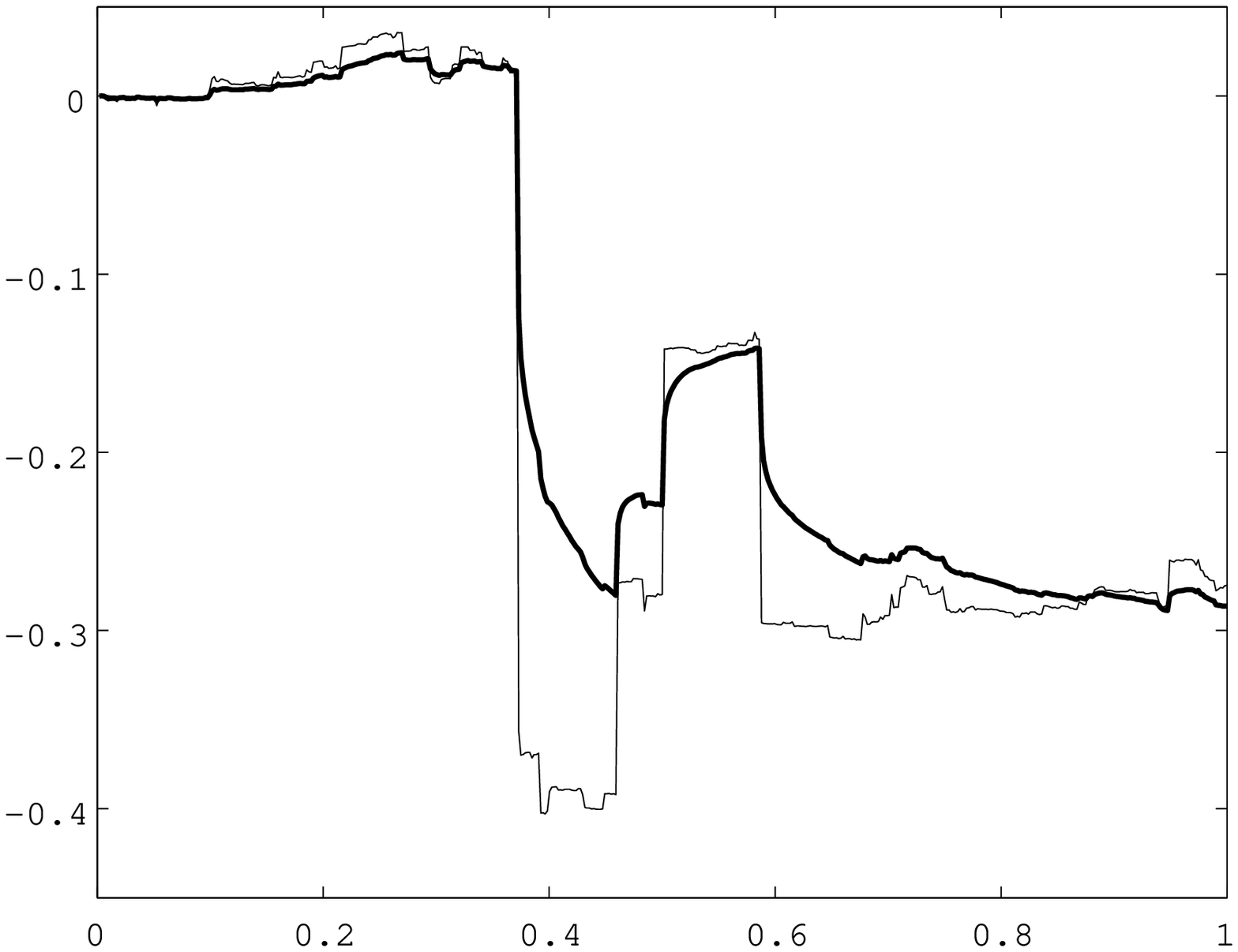}} &
      \resizebox{65mm}{!}{\includegraphics{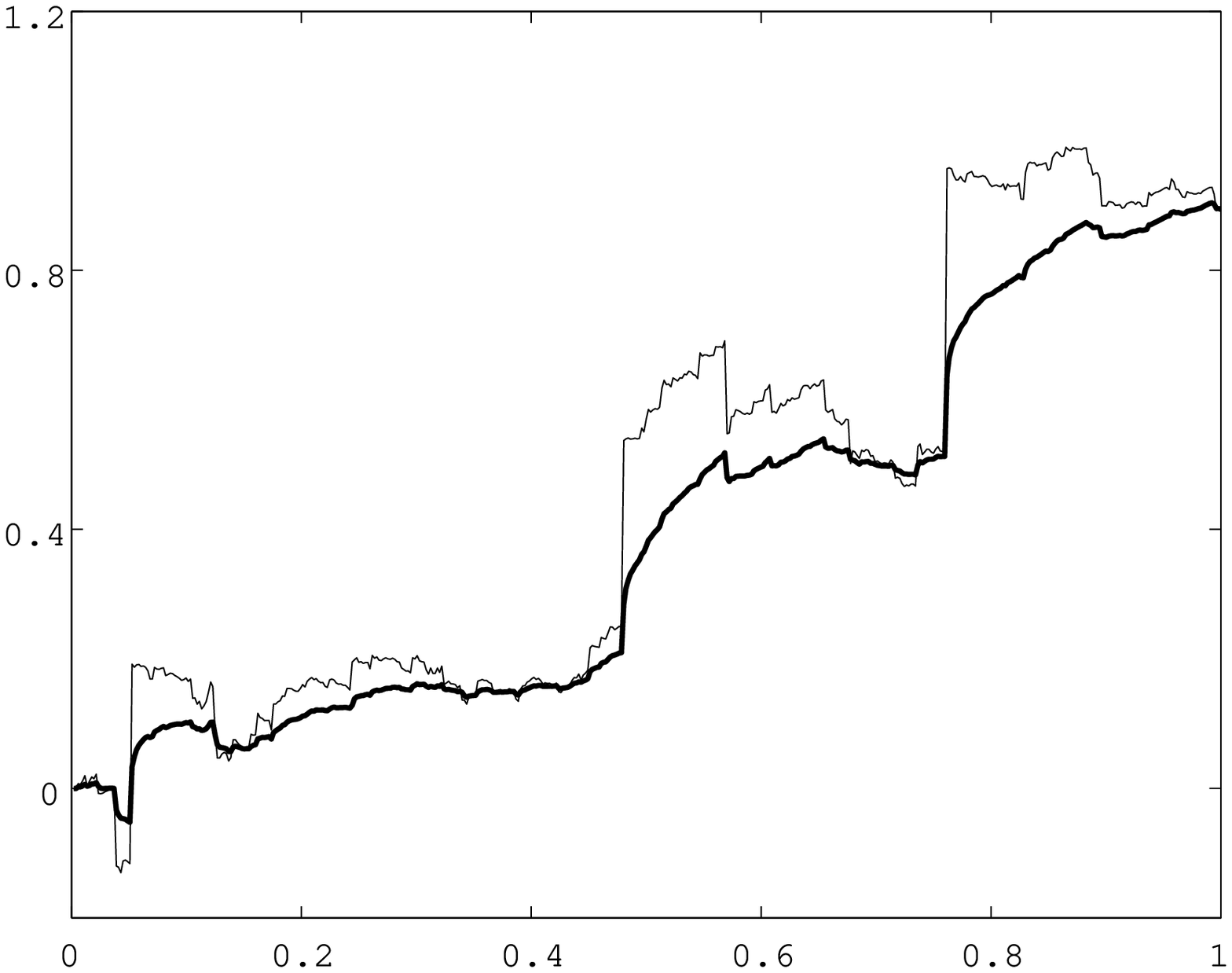}} \\
    $(H,\alpha,\rho)=(1.6,0.7,\rho_1)$ & $(H,\alpha,\rho)=(1.0,1.2,\rho_2)$ \\
      & \\
      \resizebox{65mm}{!}{\includegraphics{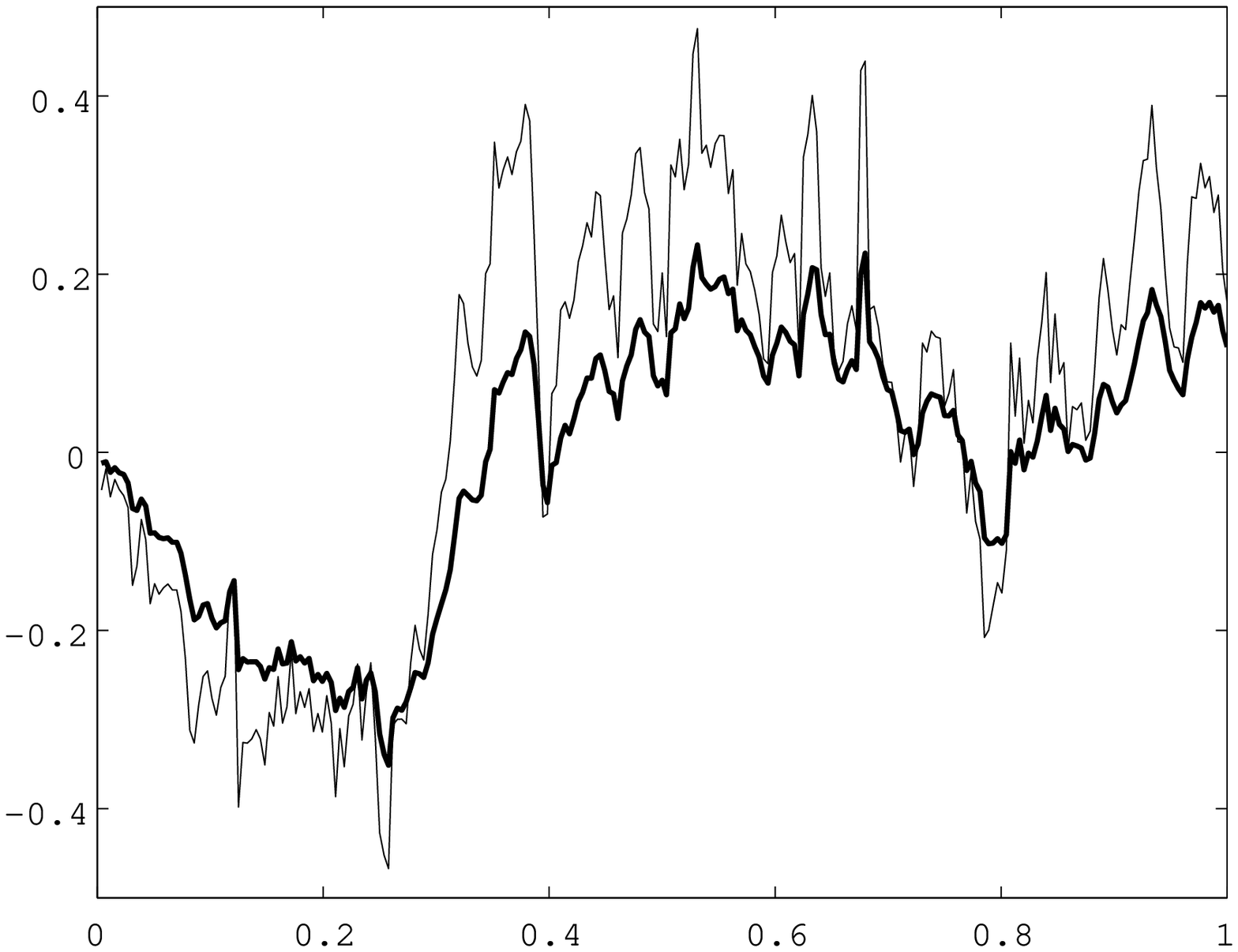}} &
      \resizebox{65mm}{!}{\includegraphics{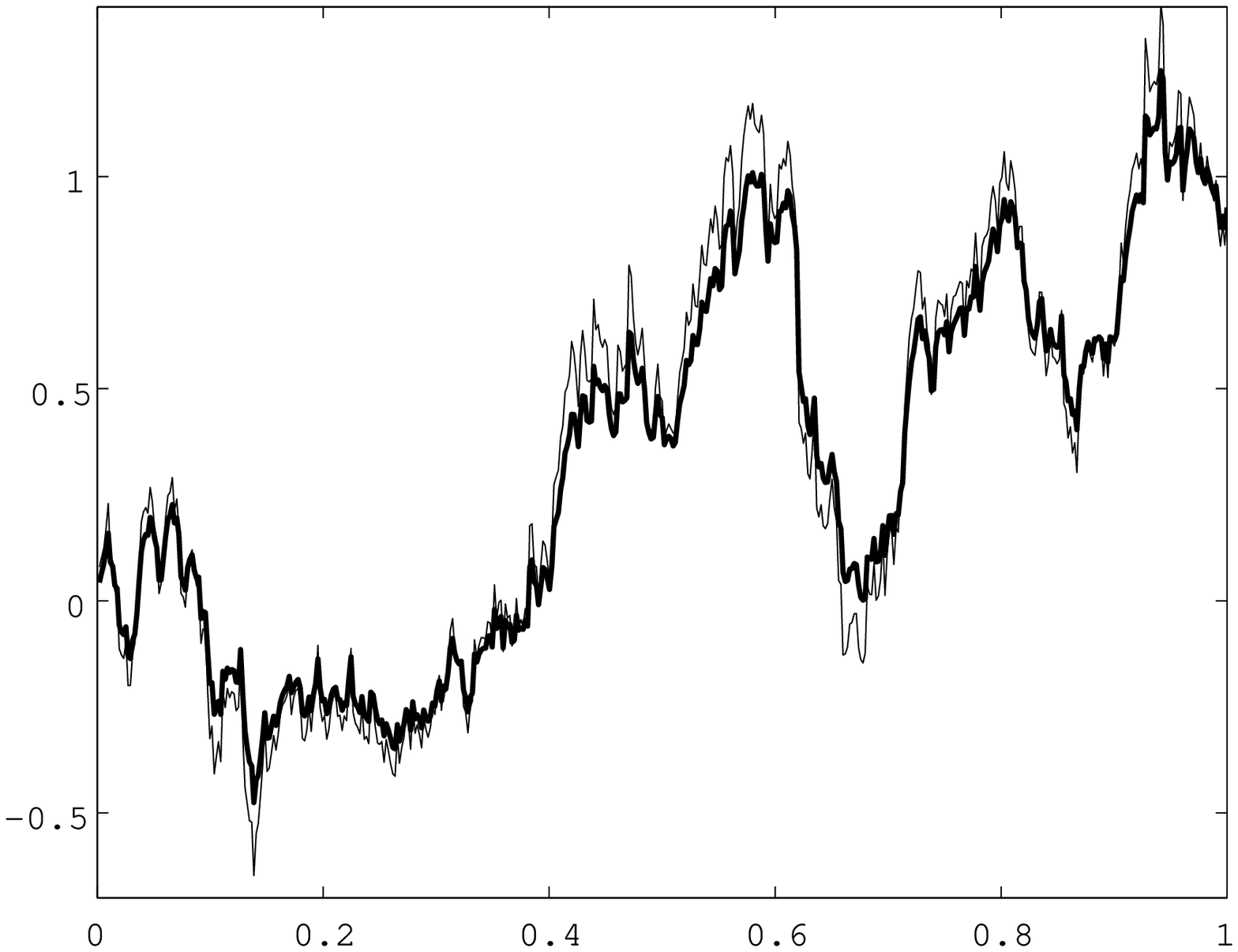}} \\
    $(H,\alpha,\rho)=(0.8,1.6,\rho_2)$ &
 $(H,\alpha,\rho)=(0.6,1.9,\rho_1)$\\
\end{tabular}
\caption{Typical sample paths of fTSm (thick line) and of its background
 driving TS process (thin line) generated via the series representation.}
\label{fTSm and TS via series representation}
\begin{tabular}{cc}
&\\
&\\
      \resizebox{60mm}{!}{\includegraphics{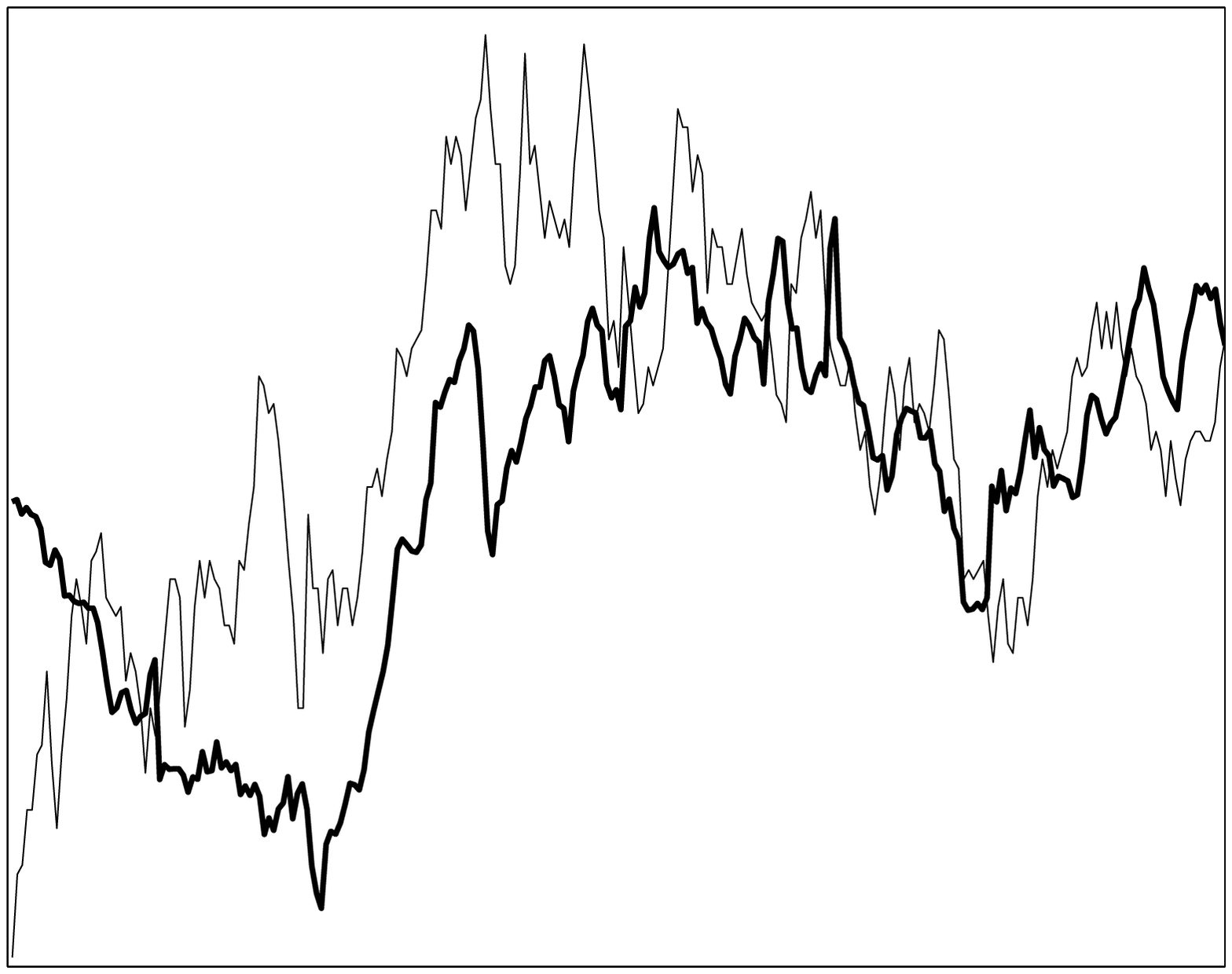}} &
      \resizebox{60mm}{!}{\includegraphics{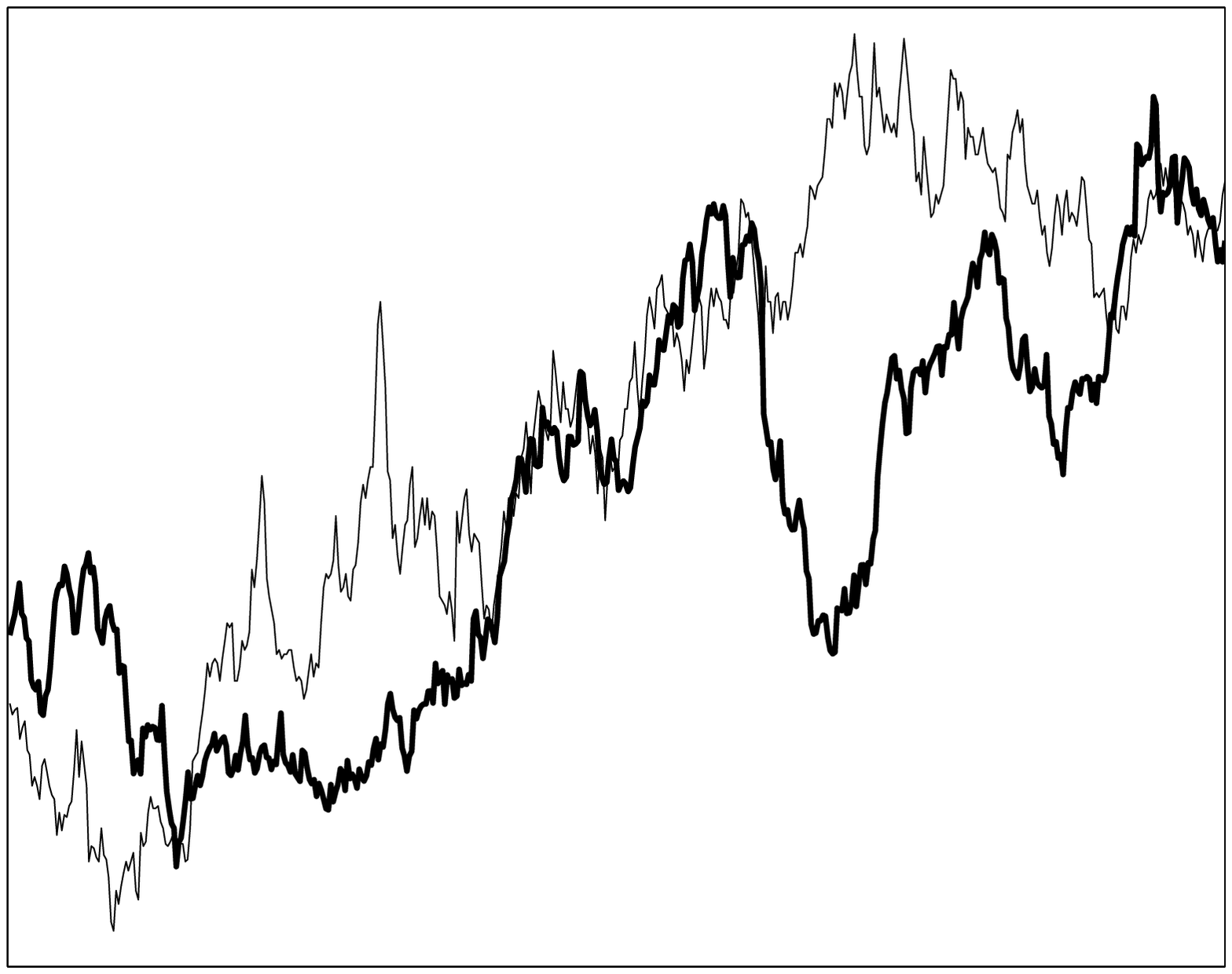}} \\
\end{tabular}
\caption{fTSm $(H,\alpha,\rho)=(0.8,1.6,\rho_2)$ (left thick)
 and $(0.6,1.9,\rho_1)$ (right thick) with (scaled) daily time series of
 TOYOTA shares on the Tokyo Stock Exchange; 247days (left thin) and
 493days (right thin) up to 02/11/2005.}
\label{TOYOTA stock price}
\end{center}
\end{figure}

In Figure \ref{fTSm and TS via series representation}, we give typical
sample paths of fTSm and of its background driving tempered stable
processes, generated via the series representation presented in
Proposition \ref{fractional Levy motion series representation}.
We put the inner measures $\rho_1$ and $\rho_2$ as
$\rho_1 (dx)=\delta_{-1.0}(dx)+\delta_{1.0}(dx)$ and 
$\rho_2 (dx)=0.5^{-\alpha}\delta_{-0.5}(dx)+\delta_{1.0}(dx).$ 
(Tempered stable L\'evy processes whose inner measure is discrete as
above are studied in Carr, Geman, Madan and Yor \cite{CGMY paper} with emphasis
on financial application and called CGMY processes.)
Observe that sample paths of fTSm look like their background driving
L\'evy process as $H$ is closer to $1/\alpha$, while the dependence
range gets longer and fTSm paths behave milder for greater $H$.

For the reader's convenience for comparison, we draw in Figure
\ref{TOYOTA stock price} daily time series of TOYOTA shares on the Tokyo
Stock Exchange, together with fTSm drawn in Figure \ref{fTSm and TS via series
representation}.
It is observed that in short time the time series looks like fTSm with $(0.8,1.6,\rho_2)$, while behaving in a Gaussian manner in longer time. 

To finish this study, let us mention that the long time behavior result
 provides yet another way to simulate sample paths of fBm.
For simplicity, consider a symmetric inner measure with a very simple
 structure, e.g. $\rho (dx)=2^{-1}(\delta_{-1}(dx)+\delta_1(dx))$, which
 reduces the random sequence $\{V_i\}_{i\ge 1}$ to a sequence of iid
 Rademacher random variables $\{\epsilon_i\}_{i\ge 1}.$ 
Observe that
\[
h^{-G}L^H_{ht}\stackrel{\mathcal{L}}{=}\sum_{i=1}^{\infty}\Bigg(m(\rho)\left(\frac{\alpha
\Gamma_i}{T}\right)^{-1/\alpha}h^{1/\alpha-1/2}\land
E_iU_i^{1/\alpha}h^{-1/2}\Bigg)\epsilon_iK_{H,\alpha}(t,T_i).
\]
Clearly, for sufficiently large $h$, the right hand side of the above
 behaves like
\[
 h^{-1/2}\sum_{i=1}^{\infty}E_iU_i^{1/\alpha}\epsilon_iK_{H,\alpha}(t,T_i).
\]
Theorem \ref{fTSm to fSm} (ii) tells us that this stochastic process (on
$[0,T]$) approximates fBm.
In order that its second moment is equal to that of a standard fBm, we set
 $h=\frac{2\alpha}{2+\alpha}N$ since then
\[
 E\left[\left(h^{-1/2}\sum_{i=1}^NE_iU_i^{1/\alpha}\epsilon_iK_{H,\alpha}
(t,T_i)\right)^2\right]=h^{-1}N\frac{2\alpha}{2+\alpha}t^{2G}=t^{2G}.
\]  
Therefore, for sufficiently large $N$, the stochastic process
\[
 \left\{\left(\frac{2\alpha}{2+\alpha}N\right)^{-1/2}\sum_{i=1}^NE_iU_i^{1/\alpha}
\epsilon_iK_{H,\alpha}(t,T_i):t\in [0,T]\right\}
\]
can be used for simulation of standard fBm.

\end{document}